\magnification=\magstep 1
\input amstex
\documentstyle{amsppt}

\def\a{\alpha}
\def\b{\beta}
\def\ga{\gamma}
\def\O{\Omega}

\def\ti{\widetilde}
\def\D{\Delta}
\def\e{\epsilon}

\def\d{\delta}
\def\A{{\Cal A}}
\def\pa{\partial}
\def\np{{\hskip 0.8pt \hbox to 0pt{/\hss}\hskip -0.8pt\partial}}
\def\G{{\Cal G}}

\def\B{{\Cal B}}
\def\F{{\Cal F}}

\def\DD{{\Cal D}}
\def\M{{\Cal M}}
\def\N{{\Cal N}}
\def\p{\phi}
\def\P{\Phi}
\def\R{{\Cal R}}
\def\si{\sigma}
\def\n{\nabla}
\def\t{\theta}
\def\lmd{\lambda}

\def\Ga{\Gamma}

\def\X{{Cal X}}

\def\GG{\bold G}

\def\NN{\bold  N}
\def\MM{\bold M}

\def\Ima{\hbox{Im}}
\def\dt{\frac{d}{dt}}

\def\SS{{\Cal S}}
\def\W{{\Cal W}}
\def\SW{\hbox{\bf SW}}
\def\sw{\hbox{\bf sw}}

\def\AustinBraamI{[3]}

\def\Chen{[5]}

\def\DonaldsonII{[8]}
\def\DonaldsonKronheimer{[9]}
\def\FloerI{[10]}

\def\FukayaI{[12]}
\def\FukayaII{[13]}
\def\GilbargTrudinger{[14]}
\def\KronheimerMrowka{[15]}

\def\Lim{[17]}
\def\Marcolli{[18]}
\def\Massey{[19]}

\def\SalamonZehnder{[21]}
\def\TaubesI{[22]}

\def\Uhlenbeck{[24]}
\def\Wang{[25]}
\def\WangYe{[26]}
\def\Witten{[27]}
\def\YeI{[28]}
\def\YeII{[29]}
\def\X{{\Cal X}}

\topmatter
\title   Equivariant and Bott-type Seiberg-Witten
Floer homology: Part I  
\endtitle
\author Guofang Wang and Rugang Ye
\endauthor
\address Max-Planck Institut f\"ur Mathematik in den Naturwissenschaften,
Inselstra{\ss}e 22-26, 04103 Leipzig, Germany
\endaddress
\email gwang\@mis.mpg.de
\endemail
\address Institute 
of System Science, Academia Sinica, 
100080 Beijing, China
\endaddress 
\email gfwang\@iss06.iss.ac.cn
\endemail
\address Department of Mathematics, University of
California, Santa Barbara, CA 93106, USA
\endaddress
\email yer\@math.ucsb.edu
\endemail
\address Ruhr-Universit\"at Bochum, Fakult\"at f\"ur Mathematik,
44780 Bochum, Germany
\endaddress
\email ye\@dgeo.ruhr-uni-bochum.de
\endemail

\rightheadtext{Bott-type  Seiberg-Witten
Floer homology  } 
\abstract
We construct Bott-type and  equivariant 
Seiberg-Witten Floer homology 
and cohomology for 3-manifolds, in particular 
rational homology spheres, and prove their
diffeomorphism invariance. This paper is a revised 
version of $\WangYe$.
\endabstract
\endtopmatter
\document

\head Table of Contents
\endhead

\noindent 1. Introduction 

\noindent 2. Preliminaries

\noindent 3. Seiberg-Witten moduli space over $Y$

\noindent 4. Seiberg-Witten trajectories: transversality

\noindent 5. Index and orientation 

\noindent 6. The temporal model and compactification

\noindent 7. Bott-type 
homology and cohomology

\noindent 8. Invariance  

\noindent Appendix A. The ordinary  Seiberg-Witten 
Floer homology

\noindent Appendix B. Local gauge 
fixing 

\noindent Appendix C. A transversality 
\bigskip

\head 1. Introduction 
\endhead

At the very beginning of the development of the  
Seiberg-Witten gauge theory it was clear that,
at least formally, the celebrated instanton homology theory of
A. Floer for 3-manifolds (homology spheres) $\FloerI$ 
could be adapted to the
Seiberg-Witten set-up. Indeed, the original
4-dimensional Seiberg-Witten equation leads naturally
to a 3-dimensional Seiberg-Witten equation via a
limit process, as first observed by Kronheimer and
Mrowka $\KronheimerMrowka$. To establish a Seiberg-Witten Floer
homology theory for a 3-manifold $Y$, the obvious idea is to replace flat 
connections in Floer's set-up by solutions of the 3-dimensional 
Seiberg-Witten equation on $Y$ (henceforth called {\it Seiberg-Witten
points}), and instanton trajectories by Seiberg-Witten
trajectories, which are solutions of the 4-dimensional Seiberg-Witten 
equation on the infinite cylinder $Y \times \Bbb 
R$. Note that the 
Seiberg-Witten points are precisely the critical 
points of the Seiberg-Witten type Chern-Simons functional, and 
that the Seiberg-Witten trajectories are 
precisely the trajectories (negative gradient flow 
lines) of this functional.  Hence the said idea amounts to establishing 
a Morse-Floer theory for the Seiberg-Witten type Chern-Simons functional.
However, one encounters various difficulties
when trying to implement this idea. The most
serious problem is that the Seiberg-Witten Floer homology for a homology
sphere (or rational homology sphere) may depend on the underlying 
Riemannian
metric, and hence is generally not a diffeomorphism
invariant. Indeed, the Euler number of the Seiberg-Witten Floer 
homology already exhibits dependence on Riemannian metrics, see 
e.g. $\DonaldsonII$ $\Chen$ and $\Lim$. The 
purpose of this series of papers, which consists 
of  the present paper (Part I), 
$\YeI$ (Part II) and $\YeII$  (Part III),
is to resolve this problem.

The trouble of non-invariance is caused by the reducible Seiberg-Witten 
point, which is 
the trivial Seiberg-Witten point in the case of a homology sphere: 
the trivial connection
coupled with the zero spinor field. It is 
a fixed point of the action of the group $S^1$ of constant 
gauges.
Under reasonable perturbations of the Seiberg-Witten equation, this
reducible Seiberg-Witten point always survives. To deal with it, one can use suitable
perturbations to make it a transversal point for
the Seiberg-Witten equation. Then one can
construct a Seiberg-Witten Floer homology, see Appendix A.
However, one encounters a serious obstruction when
trying to compare the homologies for two different perturbation parameters
(e.g. metrics). A canonical way of such comparison is to construct 
chain maps in terms of parameter-dependent Seiberg-Witten trajectories 
which connect the 3-dimensional Seiberg-Witten equation of one parameter 
to that of another. We shall call them {\it transition 
trajectories.} The said obstruction is the  presence of reducible 
transition trajectories with negative spectral flow of the
linearized Seiberg-Witten operator. Such trajectories are not in 
transversal 
position and may appear in the compactification 
of the moduli spaces of transition trajectories between 
irreducible Seiberg-Witten points. Consequently, 
the compactified moduli spaces of transition trajectories may be very 
pathological and cannot be used to define the desired chain maps. 

The appearance of such 
trajectories  can be explained in the 
following way.
The spectral flow along a reducible Seiberg-Witten
trajectory for a fixed parameter is $1$. That along a reducible 
transition trajectory from a given generic 
parameter to a nearby one is also 1. When passing from one 
generic parameter to another through certain degenerate parameters,
the spectral flow jumps and becomes negative.  
Here, typically, the Seiberg-Witten Floer homology 
also jumps. 

In many ways, the above situation is comparible with the situation of 
4-manifolds with $b^+_2=1$ where the Seiberg-Witten invariants 
exhibit dependence on chambers of Riemannian metrics, and the wall
crossing phanomena appear.  The  root of the whole trouble 
lies  in the singularity of the gauge quotient of 
the configuration space $\A(Y) \times \Ga(S)$ (see the sequel for its definition).
If we factorize the gauge group $\G$ by the subgroup 
$\G^0$ of based gauges, then 
the said singularity is seen to be the same as 
the singularity of the $S^1$ quotient 
of the space $\B^0 =(\A(Y) \times \Ga(S)) \slash \G^0$, 
where $S^1$ is the subgroup of constant 
gauges.

\bigskip 
\noindent {\bf New Constructions}
\smallskip
Since the conventional  Seiberg-Witten Floer homology
may not be diffeomorphism invariant for rational
homology spheres, we seek alternative constructions.
Our constructions are based on two ideas: one is  that 
we work on the level of the based gauge quotient, the other 
is that we multiply the based gauge quotient by a suitable 
space (e.g. the cirle), and then pass to the $S^1$ quotient. 
The first idea leads to the Bott-type theory, while the second 
leads to the equivariant theory. In the present Part I, we present 
the Bott-type theory. The equivariant theory will be presented 
in Part II.   We would like to emphasize that 
both  the Bott-type and equivariant constructions 
are  natural from the point of view of comparison with 
classical  Morse-Bott theory and equivariant homology theory.
On the other hand,  we would like to point out that as 
shown in $\Chen$ and $\Lim$ the Euler number of the ordinary Seiberg-Witten 
Floer homology can be corrected into an invariant (indeed the 
Euler number of the instanton Floer homology) by adding 
in an explicit way a certain 
(metric dependent) term to it. It is therefore 
natural to seek similar 
corrections of the entire homology.  We believe
 however that  the information
provided by
such possible corrections  is already  contained in our 
constructions. Indeed, detailed comparison 
of our constructions with the conventional construction 
should reveal the possible precise forms of desired  
corrections.
We shall discuss this point in more details in  $\YeII$.

\bigskip
\noindent {\bf The Morse-Floer-Bott Flow Complex}
 
Now, in the set-up of based gauge quotient, 
the irreducible part of the moduli space 
of gauge classes of Seiberg-Witten points consists of 
finitely many circles, while its reducible part consists of a single point, 
provided that we choose a generic parameter. These circles and the 
reducible 
point are precisely the critical submanifolds of the (Seiberg-Witten
type) Chern-Simons functional. Our goal here amounts to establishing 
a Bott-type Morse-Floer theory for the Chern-Simons
functional on the based quotient configuration space $\B^0$. 
The basic strategy is to use the moduli spaces of 
trajectories between critical submanifolds to send
(co)homological chains from one critical submanifold to 
others. Combining this map   
with the ordinary boundary operator in (co)homology 
theory then  
yields the desired boundary operator
for the (co)chain complex. This is a natural
extension of Floer's construction and was first used by Austin-Braam 
$\AustinBraamI$ and Fukaya $\FukayaI$ in Floer's set-up. The former authors use equivariant
differential forms as chains and cochains, while the 
latter uses ``geometric chains". Our core constructions 
use {\it generalized cubical singular
chains and cochains}. For delicate technical reasons 
cubical singular chains are more suitable than 
ordinary
singular chains, see $\YeI$ for details of this 
point.  For convenience, we shall sometimes 
use ``singular chains (cochains)" to refer 
to generalized cubical singular chains (cochains).

A major point of the construction using 
generalized cubical singular chains (and 
cochains) is that we can restrict to 
 subcomplexes generated by generalized 
singular cubes which are {\it $\bold F$-transversal }
(see Section 7)
with respect to  collections 
$\bold F$ of suitable maps. This is 
crucial for the invariance proof, see 
Section 8.

The compactified, suitably defined
moduli space of trajectories (flow lines) will be called the 
{\it Morse-Floer-Bott flow complex}.  The main technical 
point here  is to construct 
this flow complex along with 
its   projection to critical submanifolds,
and to establish its compactness and 
smooth structure.  In Floer's work, 
the Morse-Floer-Bott
flow 
complex does not appear. Instead, he uses the much simpler
Morse-Floer 
flow complex (indeed only the lower dimensional part of it). 
In Fukaya's work $\FukayaI$, a Morse-Floer-Bott 
flow complex appears, but 
his set-up is the full gauge quotient.
Our situation is very different, and the construction of the 
Morse-Floer-Bott flow 
complex along with its  projection to critical 
submanifolds is considerably more delicate.   We cannot find 
adquate treatments of this problem in the Seiberg-Witten
set-up or instanton set-up in the literature.

A crucial point here is how to define the endpoint 
projections of trajectories to critical submanifolds.
Obviously, one should take the limits of a trajectory 
at time infinities. The subtlety lies in the choice of 
gauges. We need to send  based gauge classes of 
trajectories to based gauge classes of critical points.
The trouble is that  based gauges on $Y \times \Bbb R$ 
may not give rise to based gauges on $Y$ at time infinities.
More precisely, if  two trajectories $u_1$ and 
$u_2$ are equivalent under based gauges, their 
limits at time infinities may not be equivalent under 
based gauges. To resolve this issue, we first transform 
a given trajectory into temporal form, and then 
take endpoint limits. This way we arrive at the 
essential 
concept of ``temporal projections". 
Two other  resulting important concepts are  
``consistent multiple temporal trajectory 
class" and 
``consistent piecewise  trajectory", which are used  
for compactifying our moduli spaces. 
The consistent condition means e.g. that the 
temporal endpoint projections of the trajectories 
in a given piecewise trajectory match each other.
We also introduce the concept ``twisted time translation".

We have two different, but 
equivalent models for our 
moduli spaces of trajectories. One is 
the {\it temporal model} $\M^0_T(S_{\a},
S_{\b})$ (see Section 6), the other is 
the {\it fixed-end model} $\M^0(p, q)$ (see Section 4).
The Fredholm theory for the Seiberg-Witten 
trajectory equation is worked out for the 
second model, and the acquired information is then 
carried over to the first one.  On the other 
hand, the temporal projections  
are   
based on the first model. 
The ordinary time translation is used for 
the temporal model, while it is  
necessary to use the twisted time translation
for the fixed-end model. 

A number of further delicate issues  have to 
be taken care of in order to compactify 
the moduli spaces and establish the smooth 
structure of the compactified moduli spaces.
 For example, in 
compactifying the moduli spaces, we have to 
establish  that after suitable gauge 
adjustments, sequences of temporal Seiberg-Witten trajectories 
of finite energy converge to consistent
piecewise 
trajectories.  
Further  discussions on the analysis in this 
regard, in particular on convergence analysis and 
gluing analysis, can be  
found in Part II $\YeI$.

\bigskip
\noindent {\bf Spinor Perturbation}
\smallskip
To prove the diffeomorphism invariance of the Bott-type
Seiberg-Witten Floer (co)homology, we employ 
the {\it transition Morse-Floer-Bott flow 
complex} (along with its projection to critical submanifolds)
 which is the flow complex built out of 
the transition trajectories.  Here, we have
to overcome the  obstruction of reducible
transition trajectories with negative spectral flow 
described before. Our strategy 
for this is to perturb the spinor equation in
the transition trajectory equation 
in order to eliminate these transition 
trajectories. We utilize the vanishing of
the rational homology group to construct suitable
vector fields which are equivariant under
based gauges. Note that {\it they are not equivariant 
under constant gauges}.  A desired perturbation is 
then gotten by adding 
one of these vector fields to 
the spinor equation. This is our key technique. Thus
 the source of our trouble,
namely the vanishing of the rational homology
group, also works to our benefit - a rather amusing 
phenomenon. 

The same kind of amusing phenomenon occurs in the 
proof of the transversality along reducible Seiberg-Witten
trajectories, see Appendix C. (Note that here we 
are not talking about transition trajectories.) This 
transversality plays a role for establishing 
the smooth structure of certain expanded moduli spaces of
transition trajectories with additional parameters,
which is used in one  stage of the 
invariance proof, see $\YeI$.

Using the transition flow complex 
with the spinor perturbation we construct 
the desired chain map from one parameter to 
another, which induces an ismorphism of the 
homologies. In this paper, we present the 
construction of this chain map, while 
leaving the  proof the isomorphism property 
to Part II and 
Part III.  Indeed, in Part II the 
isomorphism property in the equivariant 
case will be shown. In Part III, we 
shall present  an additional perturbation 
technique- the cokernel perturbation, 
and use it to prove the isomorphism 
property in the Bott-type 
case of the present paper.

We would like to mention that the 
spinor perturbations cause additional analytical difficulties which 
demand special treatments. For example, one has to establish 
a uniform $L^{\infty}$ estimate for the spinor part of the 
transition trajectories. 
With the presence of the spinor perturbations, the ordinary pointwise maximum 
principle argument no longer  works. Instead, we 
apply the 3-dimensional 
Weitzenb{\"o}ck formula (rather than the 4-dimensional one)
to obtain an initial local integral estimate in terms of 
the Seiberg-Witten energy. Then we  apply the 4-dimensional 
Weitzenb{\"o}ck formula and the technique of Moser iteration to 
derive the desired $L^{\infty}$ estimate.

\bigskip
\noindent{\bf Singular Chains/Cochains vs. Differential Forms}
\smallskip

Instead of generalized cubical singular chains and cochains on critical submanifolds, we 
can also use differential forms to build chain and cochain complexes,
and therewith a Bott-type Seiberg-Witten Floer homology and cohomology.
The boundary operator in this set-up is again constructed by
 using the Morse-Floer-Bott
flow 
complex along with its  projection to critical submanifolds. Here, the 
fibration property of the projection is essential for the construction. On the other hand, the
projection of the transition flow 
complex  to 
critical submanifolds is not a fibration in general, 
 and hence is not good enough for producing 
chain maps in the set-up of differential forms.
  In contrast, in the set-up of 
generalized cubical singular 
chains, one has the freedom of using $\bold F$-transversal
 chains for a suitable $\bold F$ and thereby bypasses 
 this problem, see Section 8.

Nevertheless, we can prove that the differential form or de 
Rham version of Bott 
theory is  diffeomorphism invariant.  Indeed, it is 
isomorphic to the singular version 
with real coefficients.  The  construction of the de Rham  version and 
the proof of its equivalence to the singular version will be given in Part II. 

\bigskip
\noindent{\bf  Uniform Formulation for All 3-Manifolds}
\smallskip

The Bott-type and equivariant constructions can easily be extended 
to 3-manifol-
ds with  nonzero first Betti number.  They can be 
shown to be isomorphic to the ordinary Seiberg-Witten 
Floer homology  
(under the assumption that the first Betti number is nonzero). 
(For accounts of 
the ordinary Seiberg-Witten Floer homolog, see $\Marcolli$,
$\Wang$ or Appendix A.)
Hence we obtain a uniform formulation of Seiberg-Witten Floer homology 
for all 3-manifolds. Details will be given in Part III. 
 In the present Part I 
and Part II, we focus on rational homology spheres.  Note that 
the ordinary Seiberg-Witten Floer homology for manifolds
 with $b_1 =1$ 
is not a diffeomorphism invariant, but its metric 
dependence is reduced 
to dependence on polarizations of the first homology.

\bigskip
A major part of the results in this paper  
were obtained in Spring  1996 while both authors were at Bochum 
University.

\head 2. Preliminaries
\endhead

To fix notations, we first recall the definitions
of the Seiberg-Witten equations on 3 and 4 dimensional manifolds. 

Let $(X_0, g)$ be an oriented Riemannian manifold of dimension
$n$ and $Spin^c(X_0)$ the set of isomorphism classses
of $spin^c$ structures on $X_0$.
Consider a $spin^c$ structure $c \in Spin^c(X_0) $ and 
its associated spinor bundle $W$ and line bundle $L$. (More
precisely, $c$ is a representative of an element in $Spin^c(X_0)$.
The homology invariants we are going to construct 
are independent of the choice of the 
representative.)
We have the associated configuration space $\A \times \Ga(W)$,
where $\A$ denotes the space of smooth unitary connections 
on $L$ and $\Ga$ the 
space of smooth sections of a vector bundle. (We suppress 
the dependence on $c $ in the notations.) The 
gauge transformation group (the group of 
gauges) is $\G = C^{\infty}(X_0, S^1)$, where $S^1 \equiv U(1)$ denotes
the unit circle in $\Bbb C$. ($\G$ depends only on $X_0$.)
The action of $\G$ on $\A \times \Ga(W)$
is defined by $((A, \P), g) \to (A + g^{-1}dg, g^{-1}\P)$.  
This formula also defines 
the (separate) actions of $\G$ on $\A$ and $\Ga(W)$. $\G$
acts freely on the subspace
of pairs $(A, \P)$ with $\P \not \equiv 0$. Such pairs
are called irreducible. The isotropy subgroup  at any reducible
pair $(A, 0)$ is the subgroup of gauge transformations 
which are constants on each component of $X_0$. If $X_0$ is 
connected, we identify it with 
$S^1$. In this case, we fix a reference point
$x_0 \in X_0$ and set $\G^0 = \{g \in \G : g(x_0)= 1\}$, which is
called the group of based gauges. Then
the quotient $\G \slash S^1$ is represented by $\G^o$. 

The action of a gauge $g$ will be denoted by 
$g^*$. We set $\B = (\A \times \Ga)  \slash \G$ and $\B^* = (\A \times
(\Ga -\{0\}))\slash \G$.  Let $\O^k(X_0)$ denote 
the space of smooth imaginary valued $k$-forms, and 
$\O^+(X_0)$ the space of smooth imaginary 
valued self-dual 2-forms (in the case that $\hbox{ dim }X_0 =4$).

We define the gauge actions on forms to be the trivial 
action, i.e. $g^* \alpha =\alpha$.  Note that the 
tangent space of the configuration space $\A \times 
\Ga(W)$ at any $(A, \P)$ can be identified with the space 
$\O^1(X_0) \times \Ga(W)$. The induced gauge action on this 
space is then the product action, more precisely, 
$g^*(\alpha, \P)= (g^*\alpha, g^* \P) = (\alpha, g^{-1} \P)$. 

We shall need the following 
\proclaim{Lemma 2.1}  
Assume that $X_0$ is closed.  Then the 
map from $\G$ to $H^1(X_0;\Bbb Z) \slash \{torsions\}$ given by $g \to 
$ the deRham class of $g^{-1}dg$ is 
surjective and induces an isomorphism 
from the component group of $\G$ to $H^1(X_0;\Bbb Z) \slash \{torsions\}$.
Moreover, there is a ~unique harmonic map $g$ with $g(x_0) = 1$ in each
component of $\G$, provided that 
$X_0$ is connected and $x_0 \in X_0$ is a fixed point.  
In particular, $\G$ is connected if $X_0$ is connected and
$H^1(X_0;\Bbb Z)$ is torsion.
\endproclaim

\demo {Proof} For simplicity,
assume that 
$X_0$ is connected.  The surjectivity of the said map
follows from integration along paths.  If $g^{-1}dg $
and $g_1^{-1}dg_1$ represent the same cohomology class,
then $g_1 = g e^{f}$ for some $f \in \O^0(X_0)$
as one easily sees. Hence $g_1$ and $g$ lie in the 
same component group. 

The statement about harmonic representative follows
from the standard theory of harmonic maps. It can 
also be derived quickly in an elementary way. For
example, if $g_1 = g e^f$ and $g$ are two 
harmonic maps, then $f$ is a harmonic function, hence
constant. \qed
\enddemo

We continue with 
the above $spin^c$ structure $c$ on $X_0$. A connection $A \in \A$ induces 
along with the Levi-Civita connection 
a connection $\n^A$ on the spinor bundle $W$ and 
the associated Dirac operator 
$D_A: \Gamma (W)
\to\Ga(W)$, 
$$D_A=\sum^n_{i=1}e_i\cdot\n^A_{e_i},$$
where $\{e_i\}$ denotes a
local orthonormal tangent frame and the 
dot denotes the Clifford multiplication. The Dirac operator is gauge 
equivariant, i.e. $D_A (g^{-1} \P) = g^{-1}D_A \P$,
and satisfies 
the following fundamental Weitzenb\"ock formula for the
Dirac operator
$$D^*_AD_A\P=-\Delta _A\P+\frac{s}{4}\P-\frac{1}{2}F_A\cdot\P,
\tag 2.1
$$
where $s$ denotes the scalar curvature of $(X,g)$ and   
$F_A$ the curvature of $A$.

Now we specify to the dimension $n = 4$. 
There is a canonical decomposition $W = W^+ \oplus W^-$ of the spinor 
bundle $W$.
The Dirac operator ~splits: $D_A : \Ga (W^+) \to \Ga (W^-)$,
$D_A: \Ga (W^-) \to \Ga (W^+)$. For a positive spinor field $\P \in
\Ga (W^+)$, the curvature $F_A$ in the above Weitzenb{\"o}ck formula
reduces to its self-dual part  $F_A^+$.

\definition{Definition 2.2}
The Seiberg-Witten equation with 
the given $spin^c$ structure $c$  is
$$
\cases
F^+_A &=\frac{1}{4}\langle e_ie_j\P,\P\rangle e^i\wedge e^j,\cr
D_A\P &=0, \cr
\endcases
\tag 2.2
$$
for $(A, \P) \in \A \times \Ga(W^+)$, where 
$\{e^i\}$ denotes the
dual of $\{e_i\} $ (a local orthonormal tangent frame).
The Seiberg-Witten operator is
$$
\hbox{\bf SW}(A, \P) = (F^+_A -\frac{1}{4}\langle e_ie_j\P,
\P\rangle e^i\wedge e^j,  D_A\P).
$$
\enddefinition

It is easy to see that the Seiberg-Witten 
operator is gauge equivariant, i.e. 
$\SW(g^*A, g^*\P) = g^*\SW(A, \P)= (F^+_A -
\frac{1}{4}\langle e_i e_j \P, \P \rangle
e^i \wedge e^j, g^{-1}D_A \P)$. 
Consequently, the Seiberg-Witten equation is gauge invariant.  

Next let $(Y,h)$ be an oriented, closed Riemannian 3-manifold with metric
$h$, and $c$ a $spin^c$ structure on $Y$. We have the associated spinor 
bundle  $S= S_c(Y)$, line bundle $L_Y = L_Y(c)$ and   
the other associated spaces: $\G(Y)$, $\A(Y)$, $\B(Y)$, 
$\B^*(Y)$ etc. {\it We shall often use $\Ga(Y)$ to denote
$\Ga(S)$.}
We set $X = Y \times \Bbb R$, which will be equipped with the product 
metric and given the  orientation $(e_1, e_2, e_3, 
\frac{\pa}{\pa t})$, where
$(e_1, e_2, e_3)$ denotes a positive local orthonormal frame on $Y$. Let 
$\pi : X \to Y$ denote the projection. 
The $spin^c$ structure $c$ induces a $spin^c$ 
structure $\pi^* c$ on $X$  with the associated line bundle $L_X =
\pi^* L_Y$ and associated spinor bundles $W^+ = \pi^* S$ 
, $W^-$.  We have the following relation between the 
Clifford multiplications on $S$ and on $W^+$:
$$v\cdot \p (y) = -(\frac{\pa }{\pa t}\cdot
 v\cdot \pi^* \p) (y, 0).$$

The associated spaces for $X$ will be indicated by 
the letter $X$.   We also 
use $\Ga^+(O)$ and $\Ga^-(O)$ to denote
$\Ga(W^+|_O)$ and $\Ga(W^-|_O)$ respectively
for $O \subset X$. 
Let $i_t: Y \to X$ denote the inclusion map which sends $y \in Y$ to
$(y, t) \in X$. A connection $A \in \A(X)$  can be written as 
$$A=a(t)+f(\cdot, t)dt,$$
where $a(t)=i_t^*(A) \in \A(Y)$ and $f \in 
C^{\infty}(X, i\Bbb R)$.
We set $\p(t)= \P(\cdot, t)=
i^*_t(\P)$ for $\P \in \Ga(W^+)$.
With these notations, we have
$$
F_A=F_a+dt \wedge\frac{\pa a}{\pa t}
+d_Yf \wedge dt,
\tag 2.3
$$
$$
2F_A^+= dt\wedge (\frac{\pa a}{\pa t}-*_YF_a
-d_Y f) -*_Y( \frac{\pa a}{\pa t}- *_YF_a -df),
\tag 2.4
$$
and
$$
D_A \P=\frac{\pa}{\pa t} \cdot (\frac{\pa \p}{\pa t} +
\np_ap + f \p).
\tag 2.5
$$ 
where $F_a$ denotes the curvature of $a$, $*_Y$
the Hodge star operator w.r.t. $h$, $d_Y$ the 
differential on $Y$, and $\np_a$ the Dirac operator 
associated with the connection $a$. 
Hence
we can rewrite (2.2) as follows
$$
\cases
\frac{\partial a}{\partial t} &=*_YF_a+d_Yf+ \langle e_i\cdot\p,\p\rangle 
e^i,
\cr
\frac{\partial \p}{\partial t} &=-\np_a\p-f\p.\cr
\endcases
\tag 2.6
$$
Here $F_a$ denotes the curvature of $a$, $*_Y$ the Hodge star operator
w.r.t. $h$, $d_Y$ the exterior
differential on $Y$, and $\np_a$ the Dirac operator associated with 
the connection $a$. 

\definition{Definition 2.3}
The Seiberg-Witten energy of $(A, \P)$ 
is
$$\eqalign{
E(A, \P)
& = \int_X(|\frac{\pa\p}{\pa t}+f\p|^2+|\np \p|^2+|\frac{\pa a}
{\pa t}-d_Yf|^2+
|*_YF_a+\langle e_i\cdot\p,\p\rangle e^i|^2).\cr
}
$$
(The volume form is omitted.) One readily shows that it is gauge 
invariant. Since $F_A = F_a + dt \wedge (\frac{\partial a}{\partial t}
-d_Yf)$, we have 
$$
E(A, \P) = \int_X (|\frac{\pa\p}{\pa t}+f\p|^2+|\np \p|^2 +
|F_A + \langle e_i\cdot\p,\p\rangle *_Ye^i|^2).
$$
\enddefinition 

Using the finite energy condition one easily derives 
from  (2.6) 
the following limiting equation for a connection 
$a\in \A(Y)$ and a spinor field $\p \in \Ga(S)$
$$
\cases
*_Y F_a+\langle e_i\cdot\p,\p\rangle e^i& =0,\cr
\np_a\p & =0.\cr
\endcases
\tag 2.7
$$

\definition{Definition 2.4} The Seiberg-Witten equation 
on $Y$ with the $spin^c$ structure $c$ is defined to be (2.7).
The Seiberg-Witten operator on $Y$ is
$$
\hbox{\bf sw}(a, \p) = (*_Y F_a+\langle e_i\cdot\p, 
\p\rangle e^i,
-\np_a\p). 
$$ 

\enddefinition

As in dimension 4, the Seiberg-Witten operator $\hbox{\bf sw}$ is 
gauge equivariant.

The following lemma is obvious.

\proclaim{Lemma 2.5} Let $u=(A, \P)= (a+fdt, \p)$
be a solution of (2.2). Denote the energy of $u$ 
on a domain $\Omega$ by $E(u, \Omega)$. Then we 
have
$$
E(u, \Omega)= 2 \int_{\Omega} |\sw(a, \p)|^2.
\tag 2.8
$$
\endproclaim

We shall need the following perturbed 
Seiberg-Witten equation

$$
\cases
*_YF_a+\langle e_i\cdot\p,\p\rangle e^i& = \nabla H(a),\cr
\np_a\p + \lambda \p & =0,\cr
\endcases
\tag 2.9
$$
where  
$\lambda$ denotes a real number  and 
$\nabla H$ the  $L^2$-gradient of a $\G(Y)$-invariant 
real valued function $H $ on $\A(Y) $. (The 
$L^2$-product is given in (2.11) below.)
We have the associated perturbed 
Seiberg-Witten operator $\hbox{\bf sw}_{\lambda, H}$.
Note that $\nabla H$ is gauge equivariant and belongs
to ${ker~ } d^*$, which are consequences of 
the gauge invariance of $H$.

The classical Chern-Simons functional plays 
a fundamental role in Floer's instanton homology 
theory. Similarly, a Chern-Simons functional associated 
with the 3-dimensional Seiberg-Witten equation will 
be important in our situation. This
functional was first used by Kronheimer and Mrowka in  their
proof of the Thom conjecture $\KronheimerMrowka$.

\definition{Definition 2.6} The Chern-Simons functional with
respect to a reference connection $a_0$ is
$$\hbox {\bf cs}(a,\p)=\frac12\int_Y(a-a_0)\wedge (F_a+F_{a_0})+
\int_Y\langle\p,\np_a\p
\rangle.
$$

Let  $\lambda$ and $H$ be as above.
The perturbed Chern-Simons functional with perturbation
$(\lambda, H)$ is

$$
\eqalign{
\hbox {\bf cs}_{\lambda,  H}(a,\p)=\frac12\int_Y(a-a_0)\wedge (F_a+
F_{a_0} )+
\int_Y\langle\p,\np_a\p\rangle  \cr
 - \lmd\int_Y\langle\p,\p\rangle 
+H(a, \p). \cr}
$$

\enddefinition

Under a gauge $g$
 the perturbed Chern-Simons functional
changes as follows:
$$\hbox {\bf cs}_{\lambda, H}(g^*(a,\p))=\hbox{\bf cs}_{(
\lambda, H)}(a,\p)+
2\pi i\int_Y c_1(L(Y))\wedge g^{-1}dg.
\tag 2.10$$
This formula implies that ${cs}_{\lambda, H}$
is invariant under the identity component of
$\G(Y)$. Hence it descends to  the quotient $\B(Y)$, 
provided that $Y$ is a rational homology sphere.

We introduce an $L^2$-product 
on $\O^1(Y) \oplus \Ga(S)$:

$$
\langle(\p_1,a_1),(\p_2,a_2)\rangle_{L^2}=\int_Y(Re
\langle\p_1,\p_2\rangle+\langle a_1,a_2\rangle)
\tag 2.11
$$
Here, $\langle a_1,a_2\rangle$ denotes the (pointwise) Hermitian product.
Easy computations lead to

\proclaim{Lemma 2.7}  The $L^2$-gradient of the perturbed 
Chern-Simons functional is given by
$$ \nabla \hbox {\bf cs}_{\lambda, H} (a, \p) = -\hbox{\bf sw}_{\lambda, H}.$$
It follows that the critical points of the 
perturbed Chern-Simons functional are precisely 
the solutions of the perturbed Seiberg-Witten 
equation.
\endproclaim

Consider a solution $(A, \P)$ of the Seiberg-Witten equation 
on the product $X$. Using a suitable 
 gauge we can transform it into temporal
form. Let's assume that it is already in
temporal form, i.e. 
$f \equiv 0$ in the formula $A = a + fdt$. 
Then Lemma 2.7 and the equation (2.6) imply
($\p(t)= \P(\cdot, t)$) 
$$
\frac{\partial }{\partial t} (a, \p) = - \nabla \hbox
{\bf cs}(a, \p).
\tag 2.12
$$

Hence 
solutions of the Seiberg-Witten equation on the product 
$X$ can be interpreted as trajectories (negative gradient 
flow lines) of the Chern-Simons functional. A similar 
formula and statement hold for solutions  of 
the perturbed Seiberg-Witten equation on $X$, 
which is  
$$
\cases
F^+_A =&\frac{1}{4}\langle e_ie_j\P,\P\rangle e^i\wedge e^j
\cr 
&+  \nabla H(a)\wedge dt+*(  \nabla H(a) \wedge dt),\cr
D_A\P =& - \lambda \frac{\pa}{\pa t} \cdot \P,\cr
\endcases
\tag 2.13
$$ 
where $\frac{d}{dt} \cdot
\P$ denotes the Clifford multiplication on $X$, $A= a + fdt 
, \p(t) = \P(\cdot, t)$ as before, and $*$ means the Hodge star
operator on $X$. The operator 
$\hbox{\bf SW}_{\lambda, H}$ is defined in an obvious way.
Obviously, we can rewrite (2.13) as follows
$$
\cases
\frac{\partial a}{\partial t} &=*_YF_a+d_Yf+
 \langle e_i\cdot\p,\p\rangle 
e^i -\nabla H(a),
\cr
\frac{\partial \p}{\partial t} &=-\np_a\p-f\p-\lambda
\p.\cr
\endcases
\tag 2.14
$$

Next we introduce the perturbed Seiberg-Witten 
energy : 

$$
\eqalign{
E_{\lambda, H}(A, \P) = 
\int_X(|\frac{\pa\p}{\pa t}+f\p|^2+|\np \p + \lambda \p|^2+|\frac{\pa a}
{\pa t}-d_Yf|^2 \cr
+|*_YF_a+\langle e_i\cdot\p,\p\rangle e^i - \nabla H(a, \p)|^2) \cr
= \int_X(|\frac{\pa\p}{\pa t}+f\p|^2+|\np \p + \lambda \p|^2
\cr
+|F_A+\langle e_i\cdot\p,\p\rangle *_Ye^i 
-*_Y \nabla H(a, \p)|^2). \cr} 
\tag 2.15
$$

Note that it is invariant under the  action of $\G(X)$.
We have an analogue of Lemma 2.5.

\proclaim{Lemma 2.8} Let $u=(A, \P)=(a+fdt, \p)$ be
a solution of (2.13). Then we have
$$
E_{\lambda, H}(u, \Omega)= 2 \int_{\Omega}|\sw_{\lambda,
H}|^2.
\tag 2.16
$$
\endproclaim

\proclaim{Lemma 2.9} Assume that $ Y$ is a rational homology sphere. 
Let $A=a+fdt$  and assume that the gauge 
equivalence class of $(a,\p)$ converges to $\a,\b\in \B(Y)$ as
$t\to-\infty,\infty$ respectively, then we have
$$\eqalign{E_{\lambda, H}(A, \P) =&
2\hbox{\bf cs}_{\lambda, H}(\a)
-2\hbox{\bf cs}_{\lambda, H}(\b)+
\int_X|D_A\P|^2\cr
&+2\int_X|
F^+_A-\frac{1}{4}\langle e_ie_j\P,\P\rangle e^i\wedge e^j
- dt \wedge \nabla H(a) -*(dt \wedge  \nabla H(a))|^2.\cr}
\tag 2.17
$$
In particular, there holds for a solution 
$(A, \P)$ of (2.13) 

$$E_{\lambda, H} (A, \P)=2\hbox{\bf cs}_{\lambda, H}(\a)
-2\hbox{\bf cs}_{\lambda, H}(\b).
\tag 2.18
$$
\endproclaim

\demo{Proof} We have
$$\int_X|D_A\P|^2
+2\int_X|
F^+_A-\frac{1}{4}\langle e_ie_j\P,\P\rangle e^i\wedge e^j
- dt \wedge \nabla H(a) -*(dt \wedge  \nabla H(a))|^2 $$ 
$$=\int_X(|\frac{\pa\p}{\pa t}+f\p+\np \p + \lambda \p|^2+|-\frac{\pa a}
{\pa t}+d_Yf+ 
*_YF_a+\langle e_i\cdot\p,\p\rangle e^i - \nabla H(a, \p)|^2)$$ 
$$=E(\P,A)+2\int\langle \frac{\pa\p}{\pa t}+f\p,\np_a\p+\lmd\p\rangle
+2\langle -\frac{\pa a}
{\pa t}+d_Yf,*_YF_a  
+\langle e_i\cdot\p,\p\rangle e^i - \nabla H(a, \p)\rangle$$
$$=2\dt\int(|\p|^2+\langle a-a_o,F_a\rangle)+2H(a)+2
\int(\langle \frac{\pa\p}{\pa t},\np_a\p\rangle+\langle\frac{\pa a}{\pa t}\p,
\p\rangle ).$$
Since the metric on $X=Y\times\Bbb R$ is the product metric,
it follows that
$$\dt\np_a\p-\np_a\dt\p= \frac{\pa a}{\pa t}\p. $$
The desired conclusion follows.\qed\enddemo

\head 3. Seiberg-Witten moduli spaces over $Y$
\endhead

We continue with the $(Y, h)$ and $c$ of the last section. 
While our theory applies 
to arbitrary closed $Y$, we 
assume for convenience that $Y$ is connected. 
Fix a reference connection $a_0$. If $L(Y)$ is a trivial bundle, we 
choose $a_0$ to be the
trivial  connection.  
We have $\A(Y)=a_0+\O^1(Y)$.
We shall use the $(l,p)$-Sobolev norms (the $L^{l, p}$-norms) 
for $l \geq 0$ and $p >0$:
$$\|u\|_{l,p}=(\sum_{0\le k\le l}\int_Y|\n^ku|^p)^{1/p}.$$
Consider the Sobolev spaces $\A_{l,p}(Y)$ and  $\Ga_{l,p}(S)$, which 
are the completions of $\A(Y)$ and $\Ga(S)$ with respect to  the 
$(l,p)$-Sobolev norm respectively. Similarly, we have  
the Sobolev spaces $\Omega^k_{l,p}(Y)$. The corresponding group of gauges is 
$\G_{l+1,p}(Y)$,
which is the completion of $\G(Y)$ with respect to the 
$(l+1, p)$-Sobolev norm.

We need to make a choice of the configuration spaces
$\A_{l, p}(Y) \times \Ga_{l, p}(S)$. We require $3p/(3-lp)
>3$ or $lp > 3$, for then all elements in $\A_{l, p}(Y) 
\times \Ga_{l,p}(S)$ are continuous, and hence the 
holonomy perturbations in the sequel can be performed. Moreover,
the corresponding gauges on the product space $X = Y
\times \Bbb R$ are continuous. In particular, if we choose 
$l=2$, then we require $p >3/2$.

\definition {Definition 3.1} We have the following spaces of
Seiberg-Witten points
$$\SS\W_{l,p} = \SS \W_{h,\lambda, H, l, p} =\{(a, \p)
\in \A_{l, p}(Y) \times \Ga_{l, p}(S): \hbox{\bf sw}_{\lambda, H}((a, \p))
=0 \}
$$
and the following moduli spaces of Seiberg-Witten points
$$\eqalign{\R_{l, p}  &=\R_{h, \lambda, H, l, p} = {\Cal SW}_{l,p}\slash
\G_{l+1, p}(Y), \cr
\R^0_{l,p} &=\R^0_{h, \lambda, H, l, p} = {\Cal SW}_{l,p} \slash \G^0_{l+1,
p}.\cr}
$$
The irreducible part of e.g. $\R_{l, p}$ will be denoted by $\R^*_{l,
p}.$
\enddefinition

We choose to work with the configuration space $\A_{2, 2}(Y)
\times \Ga_{2 ,2}(Y)$. The reason for this choice 
is that the $L^2$ spaces are more convenient to work 
with from the viewpoint of characterizing the 
cokernel of the relevant operators. 
But this  is  concerned with the situation on 
the product $X$ rather than on $Y$.  Henceforth, the subscript $l$ stands for 
$(l, 2)$, e.g. $\|u\|_l = \|u\|_{l, 2},\A_l = \A_{l, 2}.$ 
We have the quotients $\B_{2}(Y)$ and 
$\B^0_{2}(Y)$ of the chosen configuration 
space $\A_2(Y) \times \Ga_2(Y)$ under the group of 
gauges $\G_3(Y)$ and 
group of based gauges $\G_3^0$. The gauge class of 
$(a, \p)\in \A_{2}(Y) \times \Ga_{2}(S)$ with 
respect to the full gauge group $\G_3(Y)$ will be 
denoted by $[a, \p]$. 
Its gauge class with respect to based gauges
will be denoted by $[a, \p]_0$.

 Local slices of 
the gauge actions are provided by the following 
lemma, which give rise to the Banach manifold structure
for $\B_2^*$ and $\B_2^0$.

\proclaim{Lemma 3.2} 
(1) Consider an irreducible pair 
$(\bar a, \bar \p) \in \A_l(Y) \times \Ga_l(Y)$. If $u=(a, \p) \in 
\A_l(Y) \times \Ga_l(Y)$ is sufficiently close 
to $(\bar a, \bar \p)$ in $L^{l, 2}$-norm, then there is a unique 
 gauge $g= e^f \in \G_{l+1}$ such that $\tilde u
 =(\ti a, \ti \p)=g^*u$ satisfies
the gauge fixing equation
$$
d_Y^*(\tilde a- \bar a) + \hbox{Im }\langle \bar \p, \tilde
\p \rangle =0.
\tag 3.1
$$
Moreover, the $L^{l+1, 2}$-norm of  $f$ can be estimated in 
terms of the $L^{l, 2}$-norm of $u-\bar u$. 

(2)  Given pairs  $(\bar a, \bar\p), (a,\p)  \in \A_l(Y) \times
\Ga_l(Y)$, there is a unique based gauge $g= e^f \in \G_{l+1}$ 
such that $\int_Y f =0$ and $\tilde u =(\ti a,
\ti \p)= g^* u$ satisfies the gauge fixing 
equation 
$$
d_Y^*(\tilde a- \bar a)=0.
\tag 3.2
$$       
Moreover, the $L^{l+1,2}$-norm of $f$ can  be estimated in 
terms of the $L^{l,2}$-norm of $a -\bar a$.
\endproclaim

\demo{Proof} (1) The said gauge fixing equation amounts to
$$
d^*_Y d_Y f + d^*_Y (a - \bar a) + \hbox{Im} \langle
\bar \p, e^{-f} \p\rangle =0.
$$
The linearization of this equation with respect to $f$ at
$(\bar a, \bar \p)$ is the following equation
$$
d_Y^* d_Y f + f \langle \bar \p, \bar \p \rangle =0.
$$
Since $\bar \p$ is nonzero, integration shows that 
this equation has only the trivial solution. The 
desired results follow from elementary elliptic 
theory and the implicit function theorem.

(2) This part is a consequence of elementary 
theory of harmonic functions.
\qed
\enddemo

\proclaim{Lemma 3.3} (1)Each element in $\R_{2}$ or
$\R^0_{2}$ can be represented 
by a smooth pair $(a, \p)$.

(2) If $\|\n H\|_{L^{\infty}} < C$  and 
$\|\n^2 H\|_{L^{\infty}} < C$ for a constant $C$, then
 $\R_{2}$ 
and $\R^0_{2}$ are compact.
\endproclaim

\demo{Proof} We present the  proof for (2), which contains 
the argument for (1). Let $(a, \p) \in \A_2(Y)
 \times \Ga_2(S)$ be 
a solution of (2.5). Applying the 3-dimensional 
 Weitzenb-
 \"ock formula, the bound  $\|\n H\|_{
L^{\infty}}<C$ and Moser's weak maximum principle (cf. the 
proof of Lemma 8.6),  we obtain 
 $\|\p\|_{L^{\infty}}<C$ for a constant $C $. (Here and in the sequel, 
 we use the same letter $C$ to denote all constants which appear 
in a priori estimation).
Since $H^1(Y,\Bbb R)=0$, by Hodge decomposition, $a-a_0=d_Y\ga+d_Y^*\d$ for
some $\ga\in\O_{0}^0(Y)$ and $\d\in\O_{0}^2(Y)$.  
 By gauge fixing, we can assume $d_Y^*(a-a_0)=0$. Note that 
we can achieve this gauge fixing by 
a based gauge. Hence $a-a_0=d_Y^*\d$. Furthermore, we can assume
that $d_Y\d=0$. Hence we have $\D\d=F_a-F_{a_0}$. Since $\|F_a-F_{a_0}\|\le
C\|\p\|^2+C\le C$, we have $\|\d\|_2 \le C$ by elliptic estimates.
This implies
$\|a\|_1 \le C$. Applying this, the second equation of (2.5) and 
elliptic estimates, we
deduce $\|\p\|_1\le C$. Higher regularity and estimates 
follow from 
elliptic estimates and imply  the desired compactness.  
\qed
\enddemo

\proclaim{Convention 3.4} Henceforth we drop the 
subscript 
$2$ in $\SS \W_{2}$, $\R_{2} $ and $\R^0_{2}$. 
\endproclaim

Obviously, $\R = \R^0 \slash S^1$, 
where $S^1$ is the group of constant 
gauges.  We deal with $\R^*$ 
 and $\R^0$ separately.

\bigskip
\noindent {\bf The moduli space $\R^*$}
\smallskip
For a given $(a,\p) \in \A_2(Y)\times\Ga_2(Y)$,  
  let $G_Y=G_{Y,(a, \p)}
: \O_1^0(Y) \to
\O^1_0(Y) \oplus \Ga_0(S)$ be the infinitesimal gauge action operator
at $(a, \p)$, i.e.
$G_Y (f)  = (d_Yf,-f \p)$. 
Let $G_Y^* = G_{Y, (a, \p)}^*: 
\O^1_2(Y)\oplus\Ga_2(S)\to\O_1^0(Y)$ be
the formal adjoint operator
 of $G_Y$ w.r.t. the inner product (2.7). 
We have
$$G_Y^*(b,\psi)=d_Y^*b+ {Im}\langle \p,\psi\rangle.$$

There is a decomposition $\O^1_2(Y) \oplus \Ga_2(S) = 
{ker~ }G_Y^* \oplus {im~ }G_Y$. To 
be more precise, we write 
$\O^1_2(Y) \oplus \Ga_2(S) =
{ker}_2 ~G_Y^* \oplus {im}_2 ~G_Y$.
It follows that the tangent space 
$T_{[a, \p]}\B^*_2(Y)$ of $\B^*_2(Y)$ 
at $[a, \p]$ is represented by 
${ker}_2 ~G_{Y,(a,\p)}^*$ (for any representative 
$(a, \p)$ in $[a, \p]$).
Indeed, the latter gives rise to a vector  
bundle ${Ker}_2 ~G_Y^*  \to \A_2\times   (\Ga_2(S) - \{0\})$, 
whose quotient bundle ${\Cal K}er_2 ~G_Y^*$ under 
the action of $\G_3(Y)$ can 
be 
identified 
with  $T\B^*_2$.  
Next note that by 
the gauge invariance of 
the Chern-Simons functional and 
Lemma 2.7, 
we have 
$$
G_Y^* \hbox{\bf sw}_{\lambda, H} (a, \p) =0
$$
for any $(a, \p) \in \A_2(Y) \times 
\Ga_2(Y)$.  For this reason, and by the gauge equivariance,  
the operator $\hbox{\bf sw}_{\lambda, H}$ defines a section 
$[\hbox{\bf sw}_{\lambda, H}]$ of the quotient bundle
 ${ ker}_1 ~G_Y^*$,
whose fiber at $[a, \p]$ is represented 
by the kernel ${ker}_1 ~G^*_Y$ of $\G_Y^*$ in $\O_1(Y) \oplus
 \Ga_1(S).$ The moduli space $\R^*$ 
is precisely the zero locus of this section. 

In the sequel we omit the subscripts $\lambda, H$ in the 
notation $\hbox{\bf sw}_{\lambda, H}$.
Consider the operator $d\hbox{\bf sw}|_{(a, \p)} : \O^1_2(Y) 
\oplus 
\Ga_2(S) \to \O^1_1(Y)
\oplus \Ga_1(S)$, where $d \hbox{\bf sw }$  means the derivative, i.e. 
the tangent map of the operator $\hbox{\bf sw}$. By Lemma 2.6, 
it is formally self-adjoint. Assume $[a, \p]  \in \R_2 $. Then 
the gauge invariance of the equation $\hbox{\bf sw} =0$ implies 
$d \hbox{\bf sw} \circ G_Y =0$. It follows that $ G_Y^* \circ d 
\hbox{\bf sw} =0$.
Hence we obtain an operator 
$d \hbox{\bf sw}|_{{ker}_2 ~G^*_{Y, (a, \p)}}: {ker}_2
 ~G^*_{Y,(a, 
\p)} \to {ker}_1 ~G^*_{Y, (a, \p)}$ (for 
any representative $(a,\p)$ in $[a,\p]$). It is easy to see that 
this operator is Fredholm of index zero
and represents the linearization of the section 
$[\hbox{\bf sw}]$.
Let it be denoted by $\DD = \DD_{a,\p}$.
 
Lemma 2.7 implies that the operator $\DD$ coincides with the Hessian operator 
of the Chern-Simons functional with respect to the product (2.11).
Note that it extends straightforwardly to 
reducible Seiberg-Witten points. On the other hand, $\DD$ is 
gauge equivariant, i.e.
$$
d \sw |_{(g^* a, g^* \p)} = g^* d \sw |_{(a, \p)} \circ (g^{-1})^*.
$$
This is a simple consequence of the gauge equivariance of the operator 
$\sw$.  It follows that $\DD$ induces a section of Fredholm operators 
over $Y$, where the corresponding domain and target vector bundles 
are the quotient bundles induced from $\O^1_2(Y) \oplus \Ga_2(S)$
and $\O^1_1(Y) \oplus \Ga_1(S)$
respectively.  This point of view  provides a conceptual set-up.
It will appear again in the 4-dimensional discussions in the 
sequel and be used more explicitly there.

\proclaim{Definition 3.5} Let  
$(a, \p)$ be a Seiberg-Witten point, i.e.
a solution of the (perturbed) Seiberg-Witten equation (2.5). 
It is called non-degenerate if
$\DD_{a,\p}$ is onto. The classes $[a,\p]$ or $[a, \p]_0$ are 
called nondegenerate if a representative is nondegenerate.
$($This is independent of the choice of the representative.$)$
\endproclaim

\proclaim{Lemma 3.6} If all elements in $\R^*$ are nondegenerate,
then it is a naturally oriented smooth manifold of dimension zero. $($The 
orientation means that every point in $\R^*$ is assigned a sign.$)$
\endproclaim

\demo{Proof}
By the above discussions, we only need to produce the natural orientation. We can 
use either the degree of the operator $\hbox
{\bf sw}$ or the spectral flow of the operator 
 $Q$ below as in $\TaubesI$. 
(They give the same orientation.)
\enddemo

To analyse the operator $\DD$, we introduce another closely 
related formally self-adjoint Fredholm operator $Q$. 
(The Fredholm
property of $\DD$ is also a consequence of the Fredholm
property of $Q$.) First notice the following deformation complex 
$$0\longrightarrow\O^0\mathop{\longrightarrow}^{G_{a,\p}}\O^1
\oplus \Ga 
\mathop{\longrightarrow}^{d\hbox{\bf sw}}
\O^1\oplus \Ga \mathop{\longrightarrow}^{G_{a,\p}^*}
\O^0\longrightarrow 0\tag 3.3$$
where the letters $Y$ and $S$ and the Sobolev subscripts are omitted 
in the notations.
We define $Q=Q_{(a,\p)}:(\O^1_2(Y)\oplus\Ga_2(S))\oplus
\O^0_2(Y)
\to(\O^1_1(Y)\oplus\Ga_1(S))\oplus\O^0_1(Y)$  by the following 
formula:
$$Q=\pmatrix d\hbox{\bf sw} & G\cr
             G^* & 0\cr
\endpmatrix
$$

\proclaim{Lemma 3.7}  Let $(a, \p) \in \SS\W$. 
Then we have
$${ker~}Q\cong \cases\hbox{ker~}\DD\oplus\Bbb R,&
\hbox{if } (a, \p) \hbox{ is  reducible};\cr
{ker~}\DD, &\hbox{if } (a, \p)\hbox{ is irreducible}.\cr\endcases
$$
and
$${coker~}Q\cong \cases{coker~}\DD\oplus\Bbb R,&
\hbox{if } (\p,a) \hbox{ is reducible};\cr
{coker~}\DD, &\hbox{if } (\p,a)\hbox{ is irreducible}.\cr\endcases
$$
\endproclaim

We omit the simple proof. 

\bigskip
\noindent {\bf The moduli space  $\R^0$.}
\smallskip
 The above treatment does not apply to the reducible elements of $\R$, 
 because the tangent bundle
of $\B^*_2(Y)$ does not extend smoothly across the reducibles.
 To analyse the structure of $\R^0$ around reducibles, one can use 
a quotient bundle formulation on the level of the based gauge quotient. 
But we choose a different approach which gives somewhat 
stronger results. Henceforth we make

\proclaim{Assumption 3.8} $Y$ is a rational homology sphere, i.e. its 
first Betti number is zero. \endproclaim

\proclaim{Lemma 3.9} There is a canonical diffeomorphism from 
$\Sigma \equiv (a_0 + {ker}_2~d^*) \times \Ga_2(S)$ 
to $\B^0_2(Y)$.  In other words, the
former space is a global slice of the action 
of the group $\G_3^0(Y)$ on the space $\A_2(Y) \times \Ga_2(S)$.
\endproclaim 
\demo{Proof} To show that the
natural map from the former space to the latter is one to one, 
consider $b_1, b_2 \in {ker}_2 ~d_Y^*$ and $\p_1, \p_2 \in
 \Ga_2(S)$ such that $(a_0 + b_2,
\p_2) = g^*(a_0 + b_1, \p_1)$ for some gauge $g \in \G^0_3$. 
 Then $d_Y^*(g^{-1}d_Yg)=0$ and $g(y_0)=1$.  Since $Y$ is a rational 
 homology sphere, we have $dg
\equiv 0$, and hence $g \equiv 1$.  The remaining part of
 the proof is obvious.  
\qed
\enddemo
 
  This lemma enables us to reduce the Seiberg-Witten operator to the
said global slice.  But the operator $Q$ is no longer suitable for
analysing the linearization of the Seiberg-Witten operator.  Instead, 
we consider the following augmented Seiberg-Witten equation

 $$ \cases *_YF_a+d_Yf +\langle e_i\cdot\p,\p\rangle
  e^i& = \nabla H(a),\cr \np_a\p + \lambda \p +f \p & =0,\cr
  \endcases
\tag 3.4 $$
 where $a \in \A_2(Y), \p \in \Ga_2(S)$ and $ f \in \O^0_2(Y)$.

\proclaim{Lemma 3.10} Let  $(a, \p, f)$ be a solution of (3.2).
Then $(a, \p)$ satisfies the Seiberg-Witten equation (2.5) and 
$f$ is a constant. Moreover, if $(a, \p)$ is irreducible, then 
$f$ must be  zero. 
\endproclaim
\demo{Proof} Applying $d_Y^*$ to the first equation of (3.2), we 
deduce
$$
d_Y^*d_Yf+ f |\p|^2 =0.
$$
The desired conclusion follows.
\qed
\enddemo

We denote the left hand side of (3.4) by
$\hbox{\bf swa}(a, \p, f)$. The linearization
of the restriction of {\bf swa} to $\Sigma$
will be denoted by $\DD_1$. One readily checks that it is a Fredholm
operator of index 1.

\definition{Definition 3.11} Let $(a, \p)$ be a Seiberg-Witten point. 
It is called { \it based-nondege-
nerate}, if $\DD_1$ is onto at 
$(a, \p, f)$, where $f $ is an arbitrary constant 
if $\p =0$ and zero if $\p \not \equiv 0$.
It is easy to see that the based-nondegenerate property is invariant 
under gauge transformations. In particular, this 
definition makes sense for based gauge classes $[a, \p]_0$.
\enddefinition

Let the gauges act on $f$ trivially. The moduli 
space of based gauge classes of solutions of the equation (3.2) 
will be denoted by $\R^0_a$. An element in it is called
based-nondegenerate, if its corresponding element in $\R^0$
is so. As an immediate consequence of the above discussions we obtain 
\proclaim{Lemma 3.12} If all elements of $\R^0_a$ are 
based-nondegenerate, then it is a smooth oriented manifold 
of dimension one. If moreover $\R^0$ is compact, 
then the irreducible part of $\R^0_a$ consists of 
finitely many disjoint circles and its 
reducible part consists of finitely many disjoint
lines. Consequently, the irreducible part of $\R^0$
consists of finitely many disjoint circles
and its reducible part consists of finitely many points.
\endproclaim

Next we give the definition of holonomy perturbations. 
We follow $\FloerI$. Let $D$ denote 
the unit disk in $\Bbb C$. Consider a triple $(y_0, v_0, I)$,
where $y_0\in Y$, $v_0\in T_{y_0}Y$
and $I:D^2\to Y$ is a smooth embedding such that $I(0)=y_o$ and
$dI(T_0D)$ is transversal to $v_0$. 
Fix a point $s_0 \in S^1$. Let $P(y_0,v_0,I)$
be the set of all smooth embeddings  $\ga:S^1\times D\to Y$
such that $\ga(s_0,\t)=I(\t)$ for $\t\in D$ and $\frac{\pa \ga}{\pa s}
(s_0,0)=v_0$. Here $s_0$ is a fixed point in $S^1$. We set
$$P^{(m)}=\cup_{(y_o,v_0,I)}(P(y_0,v_o,I))^m,$$
for $m\in\NN$.  Now we define a map $\ga^h:\A(Y)\times
P^{(m)}\to C^{\infty}(D^2,U(1)^m)$ by
$$\ga^h(a,(\ga^1, \ga^2,\cdots,\ga^m)) (\t)=
(\ga^1_{\t}(a),\ga^2_{\t}(a),\cdot,\ga_{\t}^m(a)),$$
where $\ga^i_{\t}:\A(Y)\to U(1)$ denotes the holonomy map along 
the loop  
$\ga^i(\cdot,\t)$ (at the base point $y_0$).  It is easy to 
see that $\ga^h$ is
gauge invariant. Next  we choose a sequence $\{\e_i\}$ of 
positive numbers
as in $\FloerI$  such that
$$C^{\e}(U(1)^m,\Bbb R)=\{u\in  C^{\infty}(U(1)^m,\Bbb R)
: \| v\|_{\e}<\infty\}$$
is complete. Here 
$$\|u\|_{\e}=\sum_{i=0}^{\infty}\e_i\max_{ U(1)^m}|\n^i u|.$$

Now we set
$$\Pi=\cup_{m\in\NN}(P_m\times C^{\e}(U(1)^m,\Bbb R)).$$
This is the parameter space of 
holonomy perturbations. Choose a 
smooth function $\xi$ with support in the interior 
of $D$. For each $\pi=(\ga,u)\in \Pi$, we define
the holonomy perturbation $H_{\pi}:\A(Y)\to\Bbb R$ by
$$H_{\pi}(a)=\int_{D^2}u(\ga^h(a))\xi(\t)d^2\t.$$
It is clear that $H_{\pi}$ extends to $\A_1(Y)$.

\proclaim{Lemma 3.13} For any $\pi=(\ga,u)\in\Pi$, $H_{\pi}$ is a smooth
$\G_{3}(Y)$-invariant function. Moreover, 
the $L^2$-gradient $\n H_{\pi}$ satisfies 
$$\|\n H_{\pi}(a)\|_{L^{\infty}}\le C,$$ 
with $C>0$ independent of $a\in\A$. Similar bounds hold 
for the higher derivatives of $H_{\pi}$. The bounds can 
be made arbitrarily small by choosing $u$ small.
\endproclaim
\demo{Proof} 
For simplicity, we only consider
$m=1$. Set $H=H_{\pi}$. We can write $H(a)=\int_{D^2}u(\ga_{\t}(a))\xi 
d^2\t$. It follows that 
$$dH(a)(b)=\int_{D^2}du_{|_{\ga_{\t}(a)}}(\ga'_{\t}(a)b)\xi d^2\t.$$
Elementary computations lead to
$$\ga'_{\t}(a)b=-\ga_{\t}(a)\int_{\ga_\t}b.$$
We deduce 
$$\eqalign{dH(a)(b)&
=-\int_{D^2}\xi\langle\n u,\ga_{\t}(a)\rangle d^2\t\int_
{\ga_{\t}}b\cr
&=\int_{\ga(S^1\times D^2)}(\xi\circ\ga^{-1})\langle b,\overline{
\langle \n u,\ga_{\t}(a)\rangle}(\ga^{-1})^*(dt)\rangle |\frac{\partial\ga}
{\partial t}|^{-2} (\ga^{-1})^*(dtd^2\t)
\cr
&=\int_Y f (\xi\circ\ga^{-1})\langle b,\overline{
\langle \n u,\ga_{\t}(a)\rangle}(\ga^{-1})^*(dt)\rangle 
\cr
},$$
where
$$ f = |\frac{\partial \ga}{\partial t}|^{-2} |
(\ga^{-1})^*(dt d^2 \t) \slash dvol|.
$$ 
Consequently, $\n H(a)= f (\xi\circ\ga^{-1})\overline{
\langle \n u,\ga_{\t}(a)\rangle}(\ga^{-1})^*(dt)$. 
The desired estimate for $\n H$ follows. The higher order 
derivatives can easily be computed by using the above formula.
\qed
\enddemo

We make 
\proclaim{Assumption 3.14} 
Henceforth we choose 
$H$ in (2.5) and (3.4) to be $H_{\pi}$.
\endproclaim 

We remark in passing that for the purpose  of 
achieving transversality for the moduli 
spaces $\R^*$ and $\R^0$  it is not necessary 
to introduce the holonomy perturbations.
However, they are important for achieving transversality for 
Seiberg-Witten trajectories as will be seen in the next section.

\proclaim{Lemma 3.15} For perturbation
$\pi\in\Pi$ such that $\n^2H_{\pi}$ and 
$\n^3 H_{\pi}$ are small 
enough in $L^{\infty}$-norm (the set of such 
$\pi$ is a nonempty open set), 
there exists a unique reducible element $[(a,0)] \in\R$.
Equivalently, there is a  
unique 
$a\in\A_2(Y)$ such that
$$\eqalign{ *_YF_a-\n H_{\pi}(a) & =0,\cr
d_Y^*a & =0.\cr}
\tag 3.5$$
\endproclaim
\demo{Proof} Since $Y$ is a rational homology 
sphere, the 
operator $*_Yd_Y: {ker}~d_Y^* \to {ker}~d_Y^*$
is a bounded isomorphism. Hence the existence 
follows from the implicit function theorem. 
To prove the uniqueness,
consider connections $a$ and $a_1$ satisfying  (3.5).
We set 
$b=a-a_1$ and deduce  
$$*d_Yb=\n H_{\pi}(a_1)-\n H_{\pi}(a) \hbox{ and } d_Y^*b=0.$$
By the implicit function theorem, for $\pi$ with 
the property stated in the lemma, $b=0$.
\qed
\enddemo
The unique solution of (3.5) will be denoted by $a_{(h,\pi)}$.

\proclaim{Lemma 3.16} For $\pi\in\Pi$ satisfying the 
condition of Lemma 3.15, let
$\si(\np_{a_{(h,\pi)}})$ be the  set of eigenvalues of 
$\np_{a_{(h,\pi) }}$. 
Assume that $\n^2H_{\pi}$ and 
$\n^3 H_{\pi}$ are  small enough in $L^{\infty}$-norm
(the set of such $\pi$ is a nonempty open set). Then for 
$\lmd\in(\Bbb R-\si(\np_{a_{(h,\pi)}})) $, all 
elements of $\R(Y)$ 
are nondegenerate and all elements in $\R^0(Y)$ 
are based-nondegenerate.   
\endproclaim
\demo{Proof} 
We only present the proof for the statement 
concerning the non-degeneracy. The 
based-nondegeneracy can be treated in a similar way.
Consider $(a, \p)\in [a, \p] \in \R(Y)$, we are going to show that $\DD$ 
at 
$(a, \p)$ is onto. By gauge equivariance, we can choose $a = a_{(h, \pi)}$
for the reducible element.
By Lemma 3.7, it suffices to analyse 
the operator $Q$. 
We have 

$$Q(b, \psi, f) =\pmatrix
*_Ydb+2\langle e_i\cdot\p,\psi\rangle e^i+ d_Yf-\n^2 H(a)b\cr
\np_a\psi+\lmd\psi+f \p+ b\p\cr
d_Y^*b+{Im~}\langle\p,\psi\rangle\cr
\endpmatrix.
$$
Consider an element $(b_1, \psi_1, f_1) \in 
\O^1_1(Y) \oplus \Ga_1(S)\oplus \O_1^0(Y)$ 
satisfying
$$\langle Q(b, \psi, f), (b_1, \psi_1, f_1)\rangle_{L^2}=0\tag3.4$$
for all $(b,\psi, f)\in \O^1_2(Y)\oplus\Ga_2(S) 
\oplus \O^0_2(Y)$. 
We first derive that $(b_1,  \psi_1, f_1)$ satisfies the adjoint 
equation $Q^* =0$ (hence it satisfies $Q=0$ because $Q^*=Q$)
and is smooth.  
\smallskip
{\bf Case 1} $\p =0$ and $a= a_{(h, \pi)}$.
\smallskip
We have $\np_a \psi_1 + \lambda \psi_1 =0$. By the choice of $\lambda$, we 
conclude that $\psi_1 \equiv 0$.
Now $(b_1, f_1)$ satisfies the following equation
$$\eqalign{*_Yd_Yb_1 + d_Yf_1 - \n^2H(a)b_1 &=0, \cr
d_Y^*b_1 &= 0. \cr}
\tag3.7
$$
Since $Y$ is a rational homology sphere, the operator $(b_1, f_1) 
\to (*_Yd_Yb_1 +d_Yf_1, d_Y^*b_1$
$)$ is an isomorphism from $\O^1_2(Y) \oplus   
(\O_2^0(Y))^0$ onto $\O^1_1(Y) \oplus (\O_1^0(Y))^0$, where 
the superscript $0$ means the condition that the average be zero.
As in the proof of Lemma 3.15, we deduce that if $\nabla^2 H$
and $\nabla^3 H$ are  small enough,
$f_1$ must be a constant and $b_1 = 0$.  We conclude that 
${coker~ }Q \cong \Bbb R$. By Lemma 3.7, this implies that $\DD$
is onto.
\smallskip
{\bf Case 2} $\p \not \equiv 0$.
\smallskip
By the unique continuation, the set $U = \{\p \not = 0\}$ 
is an open dense set.  
For $y \in U$, $e_1\cdot\p(y), e_2\cdot\p(y),e_3\cdot\p(y)$ 
and $\p(y)$ span 
$S_y$, 
where
$S_y$ denotes  the fiber of $S$ at $y\in Y$. We deduce that
$\psi_1(y)=0$ for $y\in U$, whence $ \psi_1 \equiv 0$.

Now we easily see that $(b_1, f_1)$ 
satisfies the equation (3.7). Hence $b_1 \equiv 0$ and 
the equation $Q(b_1, \psi_1, f_1) =0$ reduces to 
$f_1 \p =0$.  It follows that $f_1 = 0$ in $U$ and 
consequently $f_1 \equiv 0$. We conclude that 
$Q$ is onto. By Lemma 3.7, $\DD$ is onto.
\qed
\enddemo  

As a consequence of the previous lemmas, we deduce
\proclaim{Proposition 3.17}
Let $\pi$ and $\lambda $ satisfy the same conditions as in 
Lemma 3.15 and Lemma 3.16. Then $\R$ consists of finitely many 
signed points,  the irreducible part of $\R^0$ consists of finitely 
many disjoint oriented circles, and its reducible part is a signed  point. 
\endproclaim

\definition{Definition 3.18} We call 
$\pi$ and $\lambda$ ``$Y$-generic", 
provided that they satisfy 
the conditions of Lemma 3.15 and Lemma 3.16.
The set of $Y$-generic parameters $(\pi, \lambda)$ 
is a dense open set.
\enddefinition
  
\head 4. Seiberg-Witten Trajectories: Transversality
\endhead

By Seiberg-Witten trajectories we 
mean solutions of the (perturbed) Seiberg-Witten
trajectory equation (2.13)
(we choose  $H = H_{\pi}$).  
As stated in the introduction, 
our goal is to establish a Morse-Floer 
theory 
for the Chern-Simons functional on the quotient 
space $\B^0_2(Y) = (\A_2(Y) \times \Ga_2(Y)) 
\slash  \G_3^0(Y)$.  
 The union of the critical 
submanifolds of the Chern-Simons functional is 
precisely the moduli space $\R^0$. 
The negative gradient flow
lines of the Chern-Simons functional are given by 
the temporal form of  the Seiberg-Witten trajectories, cf. Section 2.
Setting $A=a(t)+f(\cdot, t)dt$ 
and $\p(t)=\P(\cdot, t)$ in (2.10) we can rewrite it as follows
$$
\cases
\frac{\partial a}{\partial t}-*_YF_a-d_Yf-\langle e_i\cdot\p,\p\rangle e^i
&=\n H(a), \cr
\frac{\partial \p}{\partial t}+ \np_a\p+\lmd\p+f\p&=0.\cr
\endcases \tag 4.1
$$
(We omit the subscript $\pi$ in $H_{\pi}$.)

We shall use various spaces of local $(l, 2)$ Sobolev class
($L^{l, 2}_{loc}$ class), 
e.g.  $\A_{l, loc} = \A_{l, 2,loc}(X)$, $\Ga^{\pm}_{l,2, loc}$
$ =
\Ga_{l, 2,loc}(W^{\pm}) $, $\O^k_{l, loc} = \O^k_{l,2, loc}(X)$
and $\G_{l, loc}(X) = L^{l, 2}_{loc}(X, S^1)$.

\proclaim{Definition 4.1} We have the following spaces of 
Seiberg-Witten trajectories:

$\N=\{(A, \P) \in \A_{2, loc} \times \Ga_{2, loc}^+: 
(A, \P) \hbox{ solves } (4.1) $
and has finite $($perturbed$)$ 
Seiberg-Witten energy $\}$. 

The corresponding moduli spaces are:
                                   
$\M=\N
\slash \G_{3, loc}(X)$ and  
$\M^0 = \N \slash \G_{3, loc}^0(X),$  

\noindent
where $\G_{3, loc}^0(X)= \{g \in \G_{3, loc}(X): g(y_0, 0)=1\}.$
$($Recall that $y_0$ is a fixed reference point in $Y$.$)$
\endproclaim

For sets $B_1$ and $B_2$ of Seiberg-Witten points we 
set
$$\N(B_1, B_2) = \{ u \in \N: 
u(\cdot, t)
\hbox{ converges pointwise smoothly to some }  p \in B_1
\hbox{ as } t \to  $$
$$-\infty
\hbox{ and to some } q \in B_2 \hbox{ as } t \to +\infty \} 
$$
and
$$\N^G(B_1, B_2) = \{ u \in \N: \hbox{there is a }
g \in \G_{3, loc}(X) \hbox{ such that}$$ 
$$g^*u(\cdot, t)
\hbox{ converges pointwise smoothly  to some }  p \in B_1
\hbox{ as } t \to -\infty $$
$$
\hbox{and to some } q \in B_2 \hbox{ as } t \to +\infty \}. 
$$

For a positive function $\xi$ on $X$ we consider 
the following $\xi$-weighted $(l, 2)$-Sobolev norms
$$\|u\|_{l|\xi}=(\sum_{k\le l}\int_{X} \xi^2 
|\n^ku|^{2} dydt)^{\frac{1}{2}}. \tag 4.3
$$
For each  pair of nonnegative numbers
$\d = (\d_-, \d_+)$ we choose a positive smooth function
$\d_{F}$ on $\Bbb R$
such that ${\d}_F(t) =\d_{\pm}|t|$ near ${\pm} \infty$. 
We have the following $\d$-weighted $(l, 2)$-Sobolev
norms ($L^{l, 2}_{\d}$-norms)
$$
\|u\|_{l, \d}=\|u\|_{l|e^{{\d}_F}}.
$$
If $O$ is a domain in $Y \times \Bbb R$, then 
$\|u\|_{l, \d; O}$ means the $L^{l, 2}_{\d}$-norm
of $u$ on $O$. 
 
Here and in the sequel we adopt the following convention for Sobolev 
indices: in the context for $X$, the pair of indices $l, \delta$ refers to the above 
Sobolev space $L^{l,2}_{\delta}$, while
 $l, loc$ refer to the local Sobolev space $L^{l,2}_{loc}$. 
 On the other hand, in the context of $Y$, the pair of
  indices $l, p$ refers to the Sobolev space $L^{l,p}$, 
  with the index $l$ alone 
standing for the pair $l, 2$. 

 Let $\O_{l,\d}^k= \O_{l,\d}^k(X)$, 
 $\A_{l, \d}=\A_{l, \d}(X)$ and $\Ga_{l, \d}^{\pm} =\Ga_{l, \d}(W^{\pm})$ 
denote  the completion of the obvious spaces w.r.t. the $L^{l,2}_{\d}$-
norm. 

\proclaim{Proposition 4.2} Assume 
that $\lambda$ and $\pi$
are $Y$-generic. Then there are 
positive constants $\d_0$ and $E_0$ depending only
on $h$, $\lambda$ and $\pi$ 
with the following properties. For any temporal Seiberg-Witten 
trajectory $u=(A,\P)=(\p,a)$ of 
local $(2, 2)$-Sobolev class  and finite energy,
there exist a gauge $g\in \G_3(Y)$ and two smooth solutions
$u_-,u_+$ of (2.6) such that $g^*u$ is smooth
and the following holds. For all $l$,
$$\|g^*u-u_{+}\|_{l, \d_0; Y \times
[T_1, \infty)} \leq C(l),$$
$$
\|g^*u - u_-\|_{l, \d_0; Y \times [-\infty,
T_2]},
$$
where 
$C(l)$ depends only on $h$, $\lambda$, $\pi$, 
$l$ and an upper bound of the energy of $u$, and 
$T_1 >0, T_2 <0$ satisfy $E(u, Y \times [T_1, \infty)) \leq E_0,
E(u, Y \times (-\infty, T_2]) \leq E_0$. 
\endproclaim

The proof will be presented in $\YeI$.

\proclaim{Corollary 4.3} We have
$$
\M= \cup_{\a, \b \in \R} \M(\a, \b), \M^0 =\cup_{\a, \b \in \R}
\M^0(\a, \b),
$$
where $\M(\a, \b) = \N^G(\a, \b) \slash \G_{3, loc}(X)\hbox{ and }
\M^0(\a, \b)=\N^G(\a, \b)
 \slash 
\G^0_{3, loc}(X).$ 
\endproclaim

We shall use the moduli spaces $\M^0(\a, \b)$ 
(or rather suitable equivalent models for them) 
to 
construct the boundary operator in our Bott-type chain (cochain) complex.

Proposition 4.2 suggests that we can work in the set-up 
of trajectories with exponential asymptotics.

 We introduce 
the relevant 
spaces.  
For $u \in \A_2(Y) \times \Ga_2(S)$, let $G_{u}
\subset \G(Y)$ 
denote its isotropy 
 group of gauge actions. It is trivial if 
$u$ is irreducible and $S^1$ if $u$ is reducible. 
The isotropy groups are identical for gauge equivalent elements, 
hence $G_{[u]}$ is well-defined.
Choose a reference connection $a_0 \in
\A(Y)$ and set $u_0 = (a_0, 0) 
\in
\A(Y) \times \Ga(Y).$
For $p, q \in \A_2(Y) \times \Ga_2(Y)$ we 
define  the extension of $p, q$
$$
Ext(p, q)= u_0 + \chi(p-u_0) + (1-\chi)(q-u_0),
\tag 4.4
$$
where $\chi$ is a cut-off function satisfying
$$
\chi(t)= \cases 0, \hbox{ if } t \geq 0, \cr
1, \hbox{ if } t \leq -1. \cr
\endcases
\tag 4.5
$$

\definition{Definition 4.4} 
For $p, q \in \A_l(Y) \times \Ga_l(S)$ we introduce 
\smallskip
$
L_{l,\d}(p, q)= 
\{u \in \A_{l, loc} \times \Ga_{l, loc}^+: $
$u-Ext(p, q)
\in 
\O^1_{l, \d} \times \Ga^+_{l, \d}$.\}

\smallskip

$ \G_{l+1,  \d}(p, q)= \{g\in \G_{l+1, loc}(X):$ 
$g-g_0\in L^{l+1, 2}_{\d}(X, \Bbb C) 
\hbox{ for some } g_0 \in \G_{l+1,  loc}$ 
which is t-independent 
and  
belongs to $G_{p} (G_q)$ near  
$\infty  (-\infty) \},$ 

$\G_{l+1, \d}^0(p, q) = \{g \in \G_{l+1, \d}(p, q):
g(y_0, 0) =1\},$

$\G_{l+1,  \d}^I= \{g\in \G_{l+1, loc}(X):
g-1\in L^{l+1,2}_{\d}(X, \Bbb C)\}, $

$\G_{l+1, \d}^{I, 0} = \{g \in \G_{l+1, \d}^I: 
g(y_0, 0)=1\}.$

Note that  
the second and third groups act freely. The first acts freely on 
the irreducible part of $L_{l, \d}(p,q)$. 
We abbreviate e.g. $L_{2, \d}(p, q)$ to 
$L_{\d}(p, q)$, and e.g. $\G_{2+1, \d}(p, q)$ to 
$\G_{\d}(p, q)$.
We have the quotients: 

$\B_{\d} =
\B_{\d}(p, q) = L_{\d}(p, q) \slash \G_{\d}$,
 
$\B_{\d}^0 = L_{\d} \slash \G_{\d}^0, $ 

$\B_{ \d}^I = L_{\d} \slash \G_{\d}^I$,
and

$\B^{I,0}_{ \d}= L_{\d} \slash \G_{\d}^{I,0}$.

A suitable pair $(A, \P)$ induces 
the following elements (gauge classes) in these spaces respectively:
  $[(A, \P)]$, $[(A, \P)]_{0}$, $[(A,\P)]_I$ and $[(A, \P)]_{I,0}$.   
\enddefinition

Because $\G_3(Y)$ 
is connected, we obtain  equivalent (in terms of 
suitable gauges) spaces for 
Seiberg-Witten points $p', q'$
which are gauge equivalent to 
$p, q$ respectively.

\definition{Definition 4.5} For $B_1, B_2 
\subset \A_l(Y) \times \Ga_l(Y)$, we  set 
$$
L_{l, \d}(B_1, B_2)= \cup_{p \in B_1, q \in B_2}
L_{l, \d}(p, q).
$$
We abbreviate $L_{2, \d}(B_1, B_2)$ to $L_{\d}(B_1,
B_2)$,
\enddefinition

For example, we have 
 $L_{\d}(\a, \b)$ for $\a, \b \in \B_2(Y)$.
 This is a Hilbert manifold, indeed an 
affine Hilbert space bundle. We can 
represent it by the product 
$\a \times \b \times L_{\d}(p_0, 
q_0)$ for any chosen $p_0 \in
\a, q_0 \in \b$. Indeed, we can
send $u \in L_{\d}(p, q)$ to 
$(p, q, \ti u)$ with $\ti u =  u- \chi(p-p_0)
 -(1-\chi)(q-q_0)$.
This map defines 
the said Hilbert manifold structure.

The natural topology on $L_{l, \d}(\A_l(Y)
\times \Ga_l(Y), \A_l(Y) \times \Ga_l(Y))$
can be given in terms of the following distance.

\definition{Definition 4.6} For $u \in
L_{l, \d}(p_1, q_1), v \in L_{l, \d}(p_2, q_2)$ 
with $p_1, p_2, q_1, q_2 \in \A_l(Y) \times \Ga_l(Y)$,
 we set
$$
d_{l, \d}(u, v)^2=
\|p_1-
p_2\|_{l, 2}^2 + \|q_1
-q_2\|_{l, 2}^2 + 
\|(u-Ext(p_1, q_1))-
(v- Ext(p_2, q_2)) \|_{l, \d}^2.
$$
Convergence with repect to this 
distance will be called ``convergence 
in exponential Sobolev $(l, 2)$-norm,
or ``exponential convergence in 
Sobolev $(l, 2)$-norm. Smooth exponential 
convergence means exponential convergence 
in Sobolev $(l, 2)$-norm for all $l$.

These  distance and convergence concepts naturally 
descend to various quotient spaces. Indeed,
we define for equivalence classes
$\omega_1, \omega_2$
$$
d_{l, \d}(\omega_1, \omega_2)= inf\{d_{l, \d}(u, v):
u \in \omega_1, v \in \omega_2\}.
$$
\enddefinition

\definition{Definition 4.7} 
We set

$\G_{\d}^{\infty} =\{ g \in \G_{3,loc}(X): 
g-g_0 \in L^{3, 2}_{\d}(X, \Bbb C)$ for some 
$g_0 \in L^{3, 2}_{loc}$ which is $t$-independent 
near $\pm \infty \}$, 

$\G_{\d}^{\infty,0} = \{ g \in \G_{\d}^{\infty}:
g(y_0, 0)=1 \}$.
\enddefinition

\definition{Definition 4.8}
For $p, q \in \SS \W$, we set
$\N_{\d}(p,  q)
= L_{\d}(p, q) \cap \N$.
We have the moduli spaces of trajectories
$
\M_{\d}(p, q) = \N_{\d}(p, q) \slash\G_{\d}(p, q)$,
$\M_{ \d}^0(p, q) = \N_{\d}(p, q)
\slash\G_{\d}^0(p, q)$ and
$\M_{\d}^I(p, q) = \N_{\d}(p, q) \slash
\G_{\d}^I.
$
\enddefinition

\definition{Definition 4.9}
For $B_1, B_2 \subset \SS\W$, we set
$$
\N_{\d}(B_1, B_2) = \cup_{p \in B_1, q \in B_2}
\N_{\d}(p, q).
$$

For $\a, \b \in \R$, 
We 
have the moduli spaces:

$\M_{\d}(\a, \b) = \N_{\d}(\a, \b) \slash
\G^{\infty}_{\d}$,
$\M^0_{\d}(\a, \b)= \N_{\d}(\a, \b) \slash
\G^{\infty,0}_{\d}$.
\enddefinition

The irreducible part of e.g. $\M_{\d}$
will be denoted by $\M_{\d}^*$.
The spaces $\M_{\d}(p, q)$ and $\M_{\d}^0(p, q)$
are canonically isomorphic to $\M_{\d}([p], [q])$
and $\M_{\d}^0([p], [q])$ respectively. Moreover,
we have the following easy lemma.

\proclaim{Lemma 4.10} If $(\pi, \lambda)$ is $Y$-generic,
$p, q$ are smooth representatives of $\a$ and $\b$ 
respectively,
and $\d_-, \d_+$ are  
positive and do not exceed the exponent $\d_0$ in Proposition 4.3, then
there are canonical isomorphisms
$$
\M_{\d}(\a, \b) \cong \M(\a, \b), \M^0_{\d}(\a, \b)
\cong \M^0(\a, \b).
$$
We use these isomorphisms to topologize 
$\M(\a, \b)$ and $\M^0(\a, \b)$.
\endproclaim  
 
Thus, $\M(\a, \b), \M_{\d}(\a, \b)$ and $
\M_{\d}(p, q)$ can be viewed as three different models of the 
same space. The same holds for the 
spaces $\M^0(\a, \b)$ etc..
To analyse the structures of 
these moduli spaces, we focus on 
the set-ups $\B_{\d}^I(p, q)$ and $\M_{\d}^I(p, q)$.

For $p, q \in \A_2(Y) \times \Ga_2(S)$,  
consider the infinitesimal gauge action 
operator $G_X = G_{X, (A, \P)}: \O^0_{2, \d}
\to \O^1_{1, \d} \oplus \Ga_{1, \d}^+$ at a given 
$(A, \P) \in L_{\d}(p, q)$, $G_X(f) = (d f, -f \P)$. Let 
$G^*_X$ be the formal adjoint operator of $G_X$
{\it w.r.t.} the following inner product
$$
\langle(\P_1,A_1),(\P_2,A_2)\rangle_{\d}=
\int_X(Re
\langle\P_1,\P_2\rangle+\langle A_1,A_2\rangle)
e^{2\d_F}dydt.
\tag 4.6$$

There holds
$$
G_X^*(A', \P') = d^*_{\d}A' + Im \langle \P, \P' \rangle,
$$
where $d^*_{\d}= e^{-2\d_F} d^* 
e^{2\d_F}$ is the formal adjoint of $d$ with respect to 
(4.6). We have the following elementary lemma, which is 
analogous to Proposition 2a.1 in $\FloerI$.

\proclaim{Lemma 4.11} If $\d_-$ and 
$\d_+$ are positive and small enough, then 
$\O^1_{2, \d} \oplus 
\Ga_{2, \d}^+ = \hbox{im}_2 ~G_X \oplus {ker}_2 ~G_X^*$.
Consequently, the tangent space $T_{[A, \P]}\B_{\d}(p, q)$
is represented by ${ ker }~G_X^*$. 
\endproclaim
\demo{Proof} Choosing $\d_-$ and $\d_+$ positive and 
small enough, we can make sure that $G^*_XG_X
: \Omega^0_{2, \d} \to \Omega^0_{0, \d}$ is 
Fredholm.  Then the desired conclusions follow
easily.
\qed
\enddemo

Now we assume that $\d$ satisfies the condition of Lemma 4.10.
Let ${\Cal U}_{\d} \to \B_{\d}^I$ denote  
the quotient bundle of the trivial bundle $L_{\d}(p, q) \times 
(\O^+_{2, \d} \oplus \Ga_{2, \d}^-) 
\to L_{\d}(p, q)$ 
($\G_{3,\d}$ acts on $\O^+$ trivially and on $\Ga^-$ by
$g^*\Psi=g^{-1}\Psi$).
The (perturbed) 
Seiberg-Witten operator 
$\hbox{\bf SW } = \hbox{ \bf SW}_{\lambda, H}$ (cf. Section 2) induces 
a section $[\hbox{\bf SW}]$ of this bundle. If 
$p, q \in \SS\W$, then 
its zero locus is precisely the moduli 
space $\M_{\d}(p, q)$. 

The linearization of the section $[\hbox{\bf SW}]$ is 
given by the restriction of the operator $d\hbox{\bf SW}|_{(A, \P)}$ to 
${ker}_2G_X^*$, which will be denoted by $\DD_X=\DD_{X, (A, \P)}$. 
We introduce another closely related operator
$\F_{p, q}:\O^1_{2, \d}\oplus\Ga^+_{2, \d}\to
 \O_{1, \d}^+\oplus\Ga^-_{1, \d}\oplus\O^0_{1, \d}$ 

$$\F_{p, q}=\pmatrix d \hbox{\bf SW}\cr
-G^*_X\cr\endpmatrix.
$$

\proclaim{Lemma 4.12} If $p, q \in \SS\W, (A, \P) \in \N_{\d}(p, q)$,
then ${ ker~}\DD_X={ ker~}\F_{p, q}$ and
 ${ coker~}\DD_X={ coker~}\F_{p,q} .$ Consequently, $\DD_X$
is Fredholm iff $\F_{p,q}$ is Fredholm. If they are 
Fredholm, they have the same index.
\endproclaim
\demo{Proof} The kernel 
equality is clear. By the gauge invariance of the 
Seiberg-Witten equation, we have
$$d \hbox{\bf SW } \circ 
G_X =0.
\tag 4.7
$$
 Applying this and Lemma 4.11,
 we derive ${ im~ }\F_{p, q} = { im~ }
\DD_X \oplus G^*_X( { im}_2~ G_X)$. 
 But the second summand equals $\O^0_{1, \d}$. 
 Indeed, the cokernel of the operator $G^*_X G_X$ 
 is precisely $ker G_X$, which is trivial because  
 the exponential weight $\d$ implies vanishing 
 at infinity. 
\qed\enddemo

\proclaim{Lemma 4.13} Assume that $\d_- $ and $d_+$
are small enough. If $p$ $(q)$
is reducible, we assume in addition 
that $\d_-$ $(\d_+)$ is positive. Then 
$\F_{p, q}$ is Fredholm for all $p, q \in
\SS \W$.
\endproclaim
\demo{Proof} We first convert $\F_{p, q}$ into 
an equivalent form $\ti \F_{p, q}$. 
Let $\O^{k, Y}$ denote the subspace of $\O^k$ consisting of 
forms which do not contain $dt$. Then
 $\O^+_{2, \d}$ can be identified 
with $\O^{1, Y}_{2, \d}$. 
On the other hand, we can identify $\Ga^-_{2, \d}$
with $\Ga^+_{2, \d}$ by the multiplication 
with $-\pa / \pa t$.  Using these identifications
we obtain  
for a given $(A, \P) = (a + fdt, \p)$
$$\ti \F_{p,q}(b+\ti f dt, \psi)=\frac{\pa}{\pa t}
 \pmatrix b \cr \psi \cr \ti f 
\cr \endpmatrix -\pmatrix
 *_Yd_Yb+d_Y \ti f+2{Im}\langle e_i\p,\psi\rangle e^i-\n^2 H(a)\cdot b
\cr  -\np_a\psi-\lmd\psi-b\p-f\psi- \ti f\p\cr
d_Y^*b +2\d_F' \ti f+{Im}\langle\p,\psi\rangle\cr
\endpmatrix.$$
Hence $\ti \F_{p,q}-(\frac{\pa}{\pa t}-Q-
(0, 0, 2\d_F'))=(0, f, 0)$, where 
$Q$ was defined in Section 3.  Since $f$ decays exponentially, 
its multiplication is a compact operator.  
Now the limits of the operator $Q+(0, 0, 2\d_F')$ at $\pm \infty$ 
are formally self-adjoint.  Hence we can follow $\FloerI$ 
or $\SalamonZehnder$  to
show that $\dt-Q -2\d_F'$ is Fredholm. Consequently, $\F_{p,q}$ is Fredholm.
\qed\enddemo

Consider a $Y$-generic pair $(\pi_0, \lambda)$. Choose a neighborhood
$\Pi_0$ of $\pi_0$ such that $\Pi_0 \times \{\lambda\}$ 
consists of $Y$-generic pairs and the smallness conditions in 
Lemma 4.13 
is uniform for all $\pi\in \Pi_0$ (with $\lambda$ fixed).

\proclaim{Theorem 4.14} Assume that $\d_-$
and $\d_+$ are positive, satisfy the above 
smallness condition for all $\pi \in \Pi_0 $ 
and do not exceed the constant $\d_0$ in Proposition 4.2. 
Then for generic  
$\pi \in \Pi_0$ the following holds. For all $p, q \in \SS \W$, 
$[ \hbox{\bf SW}]$ is transversal to the zero section, and  hence
$\M_{\d}^I$ is a smooth manifold of dimension 
$ \hbox{ ind } \F_{p, q}$.  Consequently,
$\M_{\d}$ is a smooth manifold of dimension
$$\hbox{ ind }\F_{p, q} -
  max\{\hbox{dim }G_p, \hbox{dim }G_q\}$$ 
and $\M_{\d}^0$ is a smooth manifold of dimension 
$$\hbox{ ind }\F_{p, q}  -
 max\{\hbox{dim }G_p, \hbox{dim }G_q \}+ 1.$$ 
\endproclaim
\demo{Proof}  First assume that at least one of $p, q$ is irreducible. 
Then all elements of $L_{\d}(p, q)$ are irreducible.
 We extend the  bundle ${\Cal U}_{\d} \to$ 
$\B^I_{\d} \times \Pi_0$  in 
the  trivial way. Then $[\hbox{\bf SW}]$
gives rise to a section of the extended bundle.
By the Sard-Smale theorem, it suffices 
to show that this section is transversal to the zero section, 
which amounts to the surjectivity of the operator 
$\DD_X \oplus d_{\pi}\hbox{\bf SW } $ 
at all $(A, \P)$ which solve  the Seiberg-Witten
equation with parameter $\pi \in \Pi_0$ (and $\lambda$). 
By Lemma 4.12, the latter is equivalent to 
the surjectivity of  the operator $\F_{p, q}\oplus 
d_{\pi}\hbox{\bf SW}$, 
which in turn follows from the spinor part of the 
transversality argument in $\KronheimerMrowka$  (which is similar to the argument  in 
the proof of  Lemma 3.16) and Floer's transversality
argument in $\FloerI$  based on holonomy 
perturbations.  This establishes the statement about the space $\M^I_{\d}$.
The statements about $\M^0_{\d}$ and $\M_{\d}$
follow via the involved group actions.
For example, assume that $p$ is reducible.   We 
can choose a smooth imaginary valued function 
$f_0$ on $\Bbb R$ which is supported in the 
positive half $\Bbb R$ and equals $2\pi \sqrt{-1}$
near $+\infty$. Then the family of gauges 
$exp(t f_0), 0 \leq t <1$ represent the 
quotient $\G_{\d}(p, q) \slash \G^I_{\d}$.
Dividing by these gauges reduces the dimension 
by one. 

 If both $p$ and $q$ are reducible, then they  represent 
the same (unique) reducible element in $\R$. By 
Lemma 2.9, the energy of every Seiberg-Witten 
trajectory equals zero. By the gauge equivariance 
of the operator $\F_{p, q}$, we can use temporal gauges and assume that 
$p=q= (a_{(h, \pi)}, 0)$ and $(A, \P)  \equiv p$. The asserted 
transversality follows from Lemma C.1 in 
Appendix C. Both $\M^I_{\d}$ and $\M^0_{\d}$ are 
circles in this case.  Note that although there is 
a family of gauges $exp(t f_0)$
 as above corresponding to each end ($p$ 
or $q$), the dimension of $\M_{\d}$ is one less 
than that of the dimension of $\M^I_{\d}$,
rather than two less. This can be understood 
in the following way: suitable 
combinations of the two familties are equivalent to 
constant gauges, which do not change the constant solution
$(A, \P) \equiv p$.  

\qed
\enddemo

Note that by 
gauge equivariance, the transversality property 
is independent of the choice of the representatives
$p, q$ in their gauge classes.  Since $\R$ consists
of finitely many points, for 
generic $\pi$, the transversality property is shared 
by all $p, q$ with $[p], [q] \in \R$.

\definition{Definition 4.15} We shall say that those $\pi$ 
and the corresponding $\lambda$ as described above are {\it generic}. 
\enddefinition 

 By gauge equivariance, $d\SW, G_X$ and $G_X^*$
are equivariant as can easily be verified. 
More precisely, we have e.g.
$$
d \SW_{g^* u} = (g^{-1})^* d\SW g^*.
\tag 4.8
$$
By this and the homotopy 
invariance of Fredholm index, we also have

\proclaim{Lemma 4.16} For 
given $\a , \b \in \R$, $\hbox{ ind }\F_{p, q}$
is independent of the choice of $p \in \a,q \in \b$ and 
$(A, \P) \in L_{\d}(p, q)$. It is also independent of the choice of 
$(\d_-, \d_+)$ $($satisfying the smallness condition 
in Theorem 4.14$)$.
\endproclaim

Finally, we note  an important consequence of  transversality.

\proclaim{Lemma 4.17} 
Let $(\pi, \lambda)$ be generic, $p, q \in
\SS\W$ and $u \in \N_{\d}(p, q)$. Choose 
a reference $u_0 \in  L_{\d}(p, q)$. Then 
$d \SW_u$ has a right inverse $Q_u:\O^+_{1,\d}
\oplus \Ga^-_{1,\d} \to \O^1_{2, \d}\oplus \Ga^+
_{2, \d}$ with
$$ \|Q_u\| \leq C,$$
where $C$ depends only on $\|u-u_0\|_{2, \d}$
and $(\pi, \lambda)$. $Q_u$
is equivariant under gauge actions. In particular,
$\|Q_{g^*u}\|\leq \|g\|_{C^1}\|Q_u\|$ for $ g \in \G^{\infty}_{\d}
\cap C^1(X, S^1)$.
\endproclaim
\demo{Proof} 
Let ${\Cal O}_u$ denote the $L^2_{\d}$-orthogonal 
complement of $ker~\F_{p, q}$, where $\F_{p, q}
=\F_{p, q}|_u$. ) Then $\F_{p, q}|_{\Cal O}$
is a bounded isomorphism onto $\O^+_{1, \d}
\oplus \Ga^-_{1, \d} 
\oplus \O^0_{1, \d}$.  Let $\ti Q_u$ denote
its inverse.  By elliptic and Fredholm 
estimates, we derive $\|\ti Q\_u\| \leq
C$, where $C$ depends only on $\|u-u_0\|_{2, \d}$
and $(\pi, \lambda)$.  We set 
$Q_u= \ti Q_u|_{\O^+_{1, \d}
\oplus \Ga^-_{1, \d}}$.

The gauge equivariance of $Q_u$ follows from that 
of $d\SW$ and $G_X$.
\qed
\enddemo

\head 5. index and orientation
\endhead

Consider a $Y$-generic pair $(\pi, \lambda)$.
Let $O = O_{\pi, \lambda}$ be the unique reducible
 element in $\R$. 
For $\a\in\R$ we define 
$$\mu(\a)={ind~}\F_{p,q}-1,$$
where $p \in \a, q \in O$.

Note that $\mu $ can easily be extended to all elements of 
$\B_2(Y)$. It  
depends on $h, \pi$ and $\lambda$. Elementary computation shows 
$\mu(O) =0$. 

\proclaim{Lemma 5.1} For $p, q, r \in \SS\W $ there holds
$$\hbox{ind }\F_{p,r}=\hbox{ind }\F_{p,q}+
\hbox{ind }\F_{q,r}-dim~ G_q.$$ 
\endproclaim
\demo{Proof} 
This is similar to the corresponding index addition formula in 
Floer's theory $\FloerI$. Floer's argument can be applied directly. Another 
argument is as follows. Composing with exponential 
weight multiplication 
operators, we can transform the operators to Sobolev 
spaces without weight. Then the addition formula 
is the  consequence of a linear 
gluing argument. 
The term ${ dim~ }G_q$ arises 
because of  the ``jumping" across the kernel 
of the operator $d_Y^* + \Ima \langle \p, \rangle$
which is caused by the operator $(0, 0, \d_F')$.
\enddemo

\proclaim{Corollary 5.2} There holds 
$$
\hbox{ind }\F_{p, q} = \mu([p]) - \mu([q])+ \hbox{dim }G_q.
$$
Consequently, if $(\pi, \lambda)$ is 
a generic pair, then we have $\hbox{ dim } \M_{\d}(p, q) =
\mu([p]) -\mu([q])  - 
\hbox{dim }G_p,$
$\hbox{ dim } \M^0_{\d}(p, q) = \mu([p]) - \mu([q])+1 
 - \hbox{dim }G_p.$ 
\endproclaim
 
Next we study the orientation of the moduli spaces of 
Seiberg-Witten trajectories. 

\proclaim{Proposition 5.3}Assume that $(\pi, \lambda)$ 
is generic. Then $\M_{\d}^I(p,q), 
\M_{\d}(p, q)$ and $\M_{\d}^0(p, q)$
are orientable.  Indeed, their orientations are canonically determined 
after some choices are made, which will be 
given in the proof below. $($Consequently,
$\M_{\d}([p], [q]),$ $ \M_{\d}^0([p], 
[q]), \M([p], [q])$ and $\M^0([p], [q])$
are orientable.$)$ Moreover, the orientations 
are consistent with the gluing construction 
used in the proof below, namely the orientation 
of $\M^I_{\d}(p, q)$ is the same as the product 
orientation induced from gluing $\M^I_{\d}(p, r)$
to $\M^I_{\d}(r, q)$.
\endproclaim
\demo{Proof} 
Without loss of generality,
assume that $\d_+$ and $\d_-$ are equal and 
positive. The operator $\F_{p, q}$  defines a section 
of Fredholm operators over $L_{\d}(p, q)$. Since it is 
gauge equivariant, i.e. $\F^{g^*u}_{p, q} = g^* \F^{g}_{p, q}
(g^{-1})^*$, it induces  
a section of Fredholm 
operators over $\B^I_{\d}$, for which we use the same 
notation $\F_{p, q}$. 
Let $det(p,q)=det~\F_{p,q}$ be its determinant line bundle.

Next we 
fix a smooth irreducible pair $p_0=(a_0, \p_0)$ and
construct a new section of Fredholm operators $\tilde \F_{p_0,
p_0}$ over 
$L_{\d}(p_0, p_0)$ by interpolating  $\F_{p_0, p_0}$
along the time direction with the operator
$\F_0= \F_0^{(A, \P)}=(d^+ d^*_{\d}, D_A)$, 
which is defined  near time infinities.
More precisely, for every $u=(A, \P)$, we set
$$
\tilde \F_{p_0, p_1}u= \eta \F^u_{p, q} + (1-\eta) \F_0^u,
$$
where $\eta=\eta(t)$ is a nonnegative smooth function 
taking the value zero near infinities. By the arguments in 
the proof of 
Lemma 4.11 one readily deduces that 
$\tilde \F_{p_0, p_0}$ is Fredholm. On the other hand, it is 
gauge equivariant.  
It follows that $\tilde \F_{p_0, p_0}$ is 
Fredholm and gauge equivariant, and hence induces a section of 
Fredholm operators over $\B^I(p_0, p_0)$. Let its 
determinant line bundle be denoted by $\tilde det(p_0, p_0)$.

By the uniform structure of $\tilde \F_{p_0, p_0}$ near time 
infinities, we can deform it through a continuous 
family of Fredholm operator sections to the standard 
section $\F_0$, whose determinant 
line bundle is trivial.  
Hence $\tilde det_0(p_0, p_0)$ is 
trivial. An orientation of the vector space    $H^0_{\d}
\oplus H^1_{\d} \oplus H^+_{\d}$ (the homology of the complex 
associated with the operator $d^+ +d^*_{\d}$) then determines 
an orientation of $\tilde det_0(p_0, p_0)$, cf. $\Witten$.
We fix an orientation, i.e. a trivialization 
of $\tilde det_0(p_0, p_0)$.

For each $p \in \SS\W$, we interpolate $\F_{p, p_0}$
with $(d^+ + d^*_{\d}, D_A)$ near plus and minus  infinity to get 
$\tilde \F_{p, p_0}$ and $\tilde \F_{p_0, p}$ respectively. 

For $p, q \in \SS\W$
we construct an embedding by a simple gluing process 
$$ L_{\d}(p_0, p) \times L_{\d}(p, q) \times L_{\d}(q, p_0)
\to L_{\d}(p_0, p_0). 
$$ 
(Compare $\FukayaII$.)   
On the other hand, we choose reference elements 
$u_0 \in L_{\d}(p_0, p)$, 
$u_1 \in L_{\d}(p, q)$ and $u_2 \in 
L_{\d}(q, p_0)$. 
Then it is easy to show that $ u_0+(\hbox{ker }d^*_{\d}
\oplus \Ga_{2, \d}^+)$,$ u_1+(\hbox{ker }d^*_{\d}
\oplus \Ga_{2, \d}^+)$ and 
$ u_2   + 
(\hbox{ker } d^*_{\d} \times \Ga_{2, \d}^+)$ 
are global slices in $L_{\d}(p_0, p)$, 
$L_{\d}(p, q) $ and $L_{\d}(q, p_0)$ 
for the action of the groups $\G^I$ respectively.
Using them and the above embedding we obtain an embedding $\Theta$:
$$\B_{\d}^I(p_0, p)\times \B_{\d}^I(p, q)\times 
\B_{\d}^I(q, p_0) \to\B_1(p_0, p_0).$$
We have the projections $\pi_{0}$, 
$\pi_1$  and $\pi_2$
of the above product to its factors.  In addition let
$\pi_p$ and $\pi_q$ be its projections to $p$ and $q$ respectively. 
Now the index addition formula (Lemma 5.1) leads 
to an addition formula for the index 
bundle, which in turn implies a product 
formula for the determinant line bundle. 
We apply the last formula to the present situation to deduce 
$$\pi^*_0 \tilde det(p_0, p)\otimes \pi_1^* \det(p, q) 
\otimes \pi^*_2  \tilde det(q, p_0) \otimes 
\pi^*_pl_p  \otimes \pi^*_q l_q
\cong \tilde det(p_0, p_0)|_{\hbox{ im }\Theta},$$
where $l_p$ ($l_q$) is the dual of the 
kernel of the operator $d^* + \langle \p, \cdot \rangle$
at $p$ ($ q$).

The above isomorphism implies that $ \tilde det(p_0, p),
det(p, q)$ and $\tilde det(q, p_0)$ are trivial. Indeed, 
the above tensor 
product contains e.g. for $u_1 \in \B_{\d}^I(p,q)$
and $u_2  \in \B^I_{\d}(q, p_0)$
$$
 \pi^*_0 \tilde det(p_0, p)|_{\B^I_{\d}(p_0, p) 
\times \{u_1\} \times \{u_2\}} \otimes 
det(p, q)|_{u_1} \otimes  \tilde det(q, p_0)|_{u_2}  
\otimes l_p \otimes l_q,
$$
whose triviality clearly implies the triviality 
of $\tilde det(p_0, p)$.        

We choose an orientation for $l_O$ (note that 
the $l_p$'s are canonically equivalent
to each other for $p \in O$), and 
an orientation for each $\tilde det(p_0, p)$
and $\tilde det(q, p_0)$. Then the above isomorphism determines 
an orientation of $det(p, q)$,  which gives rise 
to an orientation of $\M^I_{\d}(p, q)$. The desired consistency 
follows from the construction. 

If both $p$ and $q$ are reducible, then 
$\M^0_{\d}$ is a circle generated by gauge actions,
and hence inherits a canonical orientation from the actions.  
Otherwise, the moduli  space $\M_{\d}$ 
is the quotient of $\M^I_{\d}$ by a free $S^1$ action 
if one of $p, q$ is reducible, and 
equals $\M^I_{\d}$ if neither is reducible. 
Hence the orientation of $\M^I_{\d}$
induces an orientation of $\M_{\d}$.
On the other hand, $\M_{\d}$ is the quotient 
of $\M^0_{\d}$ by another free $S^1$ action, 
hence we arrive at an orientation of $\M^0_{\d}$.
\qed
\enddemo

\head 6. The temporal model 
and compactification
\endhead

First we introduce 
$$
L_{T, \d}(B_1, B_2) =
\{ u \in L_{\d}(B_1, B_2):
u \hbox{ is temporal, i.e. }$$
$$  
u = (A, \P)=(a + f dt, 
\P) \hbox{ with } f \equiv 0\}.
$$
$L_{T, \d}(\a, \b)$ is a  Hilbert submanifold of 
${L}_{\d}(\a, \b)$. The  
group $\G_3^0(Y)$
acts on it freely. We obtain the 
quotient topology on $L_{T, \d}(\a, 
\b) \slash \G^0_3(Y)$. Although there is 
a natural smooth structure on this quotient 
space, we shall not use it. 
 Our 
moduli space $\M^0_T(S_{\a}, S_{\a})$ 
below is a subspace of $L_{T, \d}(\a, \b)
\slash \G^0_3(Y)$, but its smooth structure will
be  
established via its equivalent model
$\M^0(p, q), p \in \a, q \in b$.

Next we consider 
temporal Seiberg-Witten trajectories. 
 For  $ B_1, B_2 \subset \SS\W$ and 
positive number $\d$ we set 
$$
\N_{T,{\d}}(B_1, B_2) = \{u \in \N_{\d}(B_1, B_2):
u \hbox{ is temporal}\}.
$$
For $\a  \in \R$ let $S_{\a}$ denote its 
lift to $\R^0$.  For $\a, \b \in \R$ we  set
$$
\M_{T,{\d}}(S_{\a}, S_{\b}) =
\N_{T,{\d}}(S_{\a}, S_{\b}) \slash \G_3(Y), 
$$
$$
\M_{T,{\d}}^0(S_{\a}, S_{\b}) =
\N_{T,{\d}}(S_{\a}, S_{\b}) \slash \G_3^0(Y).
$$
The second space will play the major role in 
our constructions.

We shall only consider $Y$-generic parameters and fix a sufficiently small 
$\d$. Henceforth we
shall omit the subscript $\d$ in the various notations 
such as 
$\M^0_{T, \d}(S_{\a}, S_{\b})$. 

\proclaim{Lemma 6.1} 
Assume that $(\pi, \lambda)$ is generic. 
Then   
$\M_T^0(S_{\a}, S_{\b})$ (the 
``temproal model") are smooth manifolds 
which are  canonically 
diffeomorphic to $\M^0(p, q)$ (the 
``fixed-end model") for $p 
\in \a, q \in \b$.
\endproclaim
\demo{Proof} 
The space $\N(p, q)$ is 
a smooth Hilbert fibration over $\M^0(p, q)$,
and hence a Hilbert manifold. Consider 
 the  temporal transformation 
$T_G:  L(\a, \b) 
\to  L_T(\a, \b)$: for $u =(a + f  dt,\P)$, set
$$g_T(u) = e^{-\int^t_0f}
$$
and $T_G (u) = g_T(u)^*u$. 
The restriction of $T_G$ to $\N(p, q)$ 
has image $\N_T(\a, \b) \subset  L_T(\a, \b)$. 
The induced temporal transformation 
$T_G: \M^0(p, q) \to  \M^0_T(S_{\a}, 
S_{\b}) \subset  L_T(\a, \b) \slash 
\G^0_3(Y)$ is a homeomorphism. We use this 
natural homeomorphism to define the 
canonical smooth structure for $\M^0_T(S_{\a},
S_{\b}).$
\qed
\enddemo

\definition{Definition 6.2} For $u \in \N(p, q)$,
let $p', q'$ be the endpoints of $T_G(u)$, i.e. 
its limits  at
$-\infty$ and $+\infty$ respectively. We set
$\pi_-(u) =[p']_0, \pi_+(u) = [q']_0$. $\pi_+$
and $\pi_-$ are called the {\it temporal projections}.
For $[u]^T_0 \in \M_T^0(S_{\a}, S_{\b})$ we define 
the temporal projections as follows:
$\pi_{\mp}([u]_0^T) = \pi_{\mp}(u)$. This is 
independent of the choice of $u$. We have a similar 
definition for $\M^0_{\d}(p, q)$.
\enddefinition

\proclaim{Lemma 6.3} The temporal projections 
are invariant under the action of $\G^{\infty, 0}_{\d}$.
\endproclaim
\demo{Proof} Consider $u_1, u_2 
\in \N(p, q)$ such that $u_2 = g^* u_1$ with 
$g \in \G^{\infty, 0}_{\d}$. Then $T_G(u_1)
= \tilde g^* T_G(u_2)$ for some $\tilde 
g \in \G^{\infty, 0}_{\d}$. Since both $T_G(u_1)$
and $T_G(u_2)$ are temporal, if follows that 
$\tilde g \in \G^0_3(Y)$. Hence we conclude 
that the temporal projections of $u_1$ are 
the same as those of $u_2$.
\enddemo

We need to compactify our moduli spaces of trajectories.
For this purpose, we introduce the following definitions.

\definition{Definition 6.4} Let $\M_T^0(S_{\a_0},
..., S_{\a_k})$ denote the fibered product $$\M_T^0(S_{\a_0},
S_{\a_1}) \times_{S_{a_1}} \M_T^0(S_{a_1}, S_{\a_2}) ...
\times_{S_{\a_{k-1}}} \M_T^0(S_{\a_{k-1}}, S_{\a_k}),$$
i.e.
$$
\M_T^0(S_{\a_0}, ..., S_{\a_k})= \{ ([u_1]_0^T,...,[u_k]_0^T)
\in \M_T^0(S_{\a_0}, S_{\a_1}) \times ...
\M_T^0(S_{\a_{k-1}}, S_{\a_k}):$$
$$ \pi_+([u_i]^T_0 )=
\pi_-([u_{i+1}]_0^T), i=1, ..., k-1\}.
$$
We call elements in $\M_T^0(S_{\a_i}, S_{\a_{i+1}})$
{\it temporal trajectory classes}, elements 
in $$\M_T^0(S_{\a_0}, S_{\a_1}) \times... 
\times \M_T^0(S_{\a_{k-1}}, S_{\a_k})$$ 
{\it multiple  temporal trajectory classes}, 
and elements in $$\M^0_T(S_{\a_0}, ...,
S_{\a_k})$$  {\it consistent multiple   
temporal trajectory classes}. We call $S_{\a_i}$ the 
$i$-th {\it juncture manifold} of 
these consistent mulitple temporal trajectory 
classes.   We call $[u_i]_0^T$
the {\it $i$-th piece} of $([u_1]_0^T, ..., [u_k]_0^T)$.

The temporal projections $\pi_{\pm}$ are naturally 
extended to consistent multiple temporal trajectory 
classes. Indeed, we define $\pi_-$ in terms of 
the first piece, and $\pi_+$ in terms of 
the last piece.

For distinct $\a$ and $\b$, we define 
$
\M^0_T(S_{\a}, S_{\b})_k
$
to be the union of all $\M_T^0(S_{\a_0}, ..., $
$S_{\a_k})$
with distinct $\a_0,..., \a_k$ and $\a_0=\a, \a_k=\b$.
Finally, we set for distinct $\a, \b$
$$
\MM_T^0(S_{\a}, S_{\b})= \cup_{k} \M^0_T(S_{\a},
S_{\b})_k.
$$
By the definition, $\MM_T^0(S_{\a}, S_{\b})$
has a natural stratified structure.

For $\a =\b$, we set $\MM_T^0(S_{\a}, S_{\a})=
\M_T^0(S_{\a}, S_{\a})$, which consist of 
time-independent trajectories. 
\enddefinition

\definition{Definition 6.5}  Let $p, q \in {\Cal S}{\Cal W}$.
A  {\it  piecewise trajectory} $u = (u_1, ..., u_m) 
$ of length $k$ from $p$ to $q$ with consecutive 
{\it junctures} $p_0=p, ..., p_k=q \in \SS\W$ is an element in 
$\N(p_0, p_1) \times ...\N(p_{k-1}, p_k)$ with
 $u_i \in \N(p_{i-1}, p_i)$.
$u_i$ is called the {\it i-th piece} of $u$. 
$p, q$ are called the {\it endpoints} of $u$ at 
$+\infty$ and $-\infty$ respectively.  $u$ is called 
{\it proper}, if $E(u_i) >0$ for 
every $i$. 
We also call a piecewise trajectory of length $k$ a 
{\it $k$-trajectory}.
\enddefinition

\definition{Definition 6.6} A piecewise trajectory 
$u=(u_i, ..., u_k)$ is called {\it consistent}, 
provided that $\pi_+(u_i)= \pi_-(u_{i+1}),
i=1,..., k-1.$

For 
$p_0, ..., p_k \in \SS\W$, let  $\N(p_0, ..., p_k)$
denote the space of consistent 
$k$-trajectories  with junctures $p_0, ..., p_k$. 
For $p, q \in \SS\W$ with $[p] \not =[q]$, let  
$ \N(p, q)_k $ denote the space of proper and consistent
$k$-trajectories from $p$ to $q$. We set 
$$
{\bold N}(p, q)= \cup_k \N(p, q)_k.
$$ 
 
The temporal projections $\pi_+$ and $\pi_-$ naturally extend to 
consistent piecewise trajectories, namely $\pi_-$ is 
defined in terms of the first portion and $\pi_+$ in 
terms of the last portion.  
\enddefinition

\definition{Definition 6.7}
For  $p_0, .., p_k \in \SS\W$, we set 
$$ \M^0(p_0,..., p_k)= \N(p_0,..., p_k) \slash 
(\G_{\d}^0(p_{0}, p_1) \times...\times 
\G_{\d}^0(p_{k-1}, p_k)).$$

For each $\a \in \R$, we choose an element $p_{\a} 
\in \a$.  We fix this choice henceforth and denote 
the set of these elements by $\SS\W_0$. 
For distinct $p, q \in \SS\W_0$, 
let $ \M^0(p, q; \SS\W_0)_k$ denote the union of 
all $\M^0(p_0,..., p_k)$ with distinct $p_0=p,..., p_k=q \in
 \SS\W_0$, 
and let $\MM^0(p, q; \SS\W_0)$ denote 
the union of $\M^0(p, q;
\SS\W_0)_k$
over all possible $k$. By the definition, $\MM^0(p, q;
\SS\W_0)$
has a natural stratification.
\enddefinition

\proclaim{Proposition 6.8} Assume that 
$(\pi, \lambda)$ is generic. Consider 
$\a_i \in \R$ with $\mu(\a_i) \geq  \mu(\a_{i+1})$, 
$i =0, ..., k-1$. Also consider $p_i \in 
\a_i, i=0, ..., k$.  The moduli spaces 
 $\M^0(p_0,...,p_k)$ and $\M_T^0(S_{\a_0},
...,S_{\a_k})$ are smooth manifolds of 
dimension 
$$
\sum_{0 \leq i \leq k-1} dim~\M^0(p_i, p_{i+1})
-\sum_{1 \leq i \leq k-1} dim~ S_{\a_i}= $$
$$
 \sum_{0 \leq i \leq k-1} dim~\M^0_T(S_{\a_i}, S_{\a_{i+1}})
-\sum_{1 \leq i \leq k-1} dim~ S_{\a_i}. 
$$
Moreover, $\M^0(p_0, ..., p_k)$ is 
canonically diffeomorphic to   $\M_T^0(S_{\a_0},
...,S_{\a_k})$. 
\endproclaim
\demo{Proof} Consider e.g. the case $\M^0_T(S_{\a_0}, 
S_{\a_1}, S_{\a_2})$. The general case is similar.
We have the temporal  projections 
$\pi_+: \M^0_T(S_{\a_0}, S_{\a_1}) \to S_{\a_1}$ and 
$\pi_-: \M^0(S_{\a_1}, S_{\a_2}) \to S_{\a_1}$. They  
are smooth maps and transversal to each 
other because of the action of the constant gauges.
The moduli space $\M_T^0(S_{\a_0},
S_{\a_1}, S_{\a_2})$ is precisely the 
preimage of the digonal in $S_{\a_1} \times
 S_{\a_1}$ under the map $ \pi_+ \times 
 \pi_-$. Hence we infer that it is a smooth 
 manifold. The dimension formula follows easily.
 
 The proof for $\M^0(p_0, ..., p_k)$ is 
 similar. Here constant gauges are replaced by
 the family of gauges $e^{t f_0}, 0\leq t < 1$, where 
 $f_0$ is a compactly supported imaginery 
 valued function on $\Bbb R$ with $f_0(y_0, 0)=
 2\pi \sqrt{-1}$. 
 
 The claimed diffeomorphism 
 is given by the temporal transformation.
 \qed
 \enddemo

Next we introduce suitable convergence concepts.
First, we set for $p_i \in 
\A_l(Y) \times \Ga_l(Y)$,
$$
 L_{l, \d}(p_0, ..., p_k) = L_{l, \d}(p_0, p_1) 
 \times L_{l, \d}(p_1, p_2) \times ... L_{l, \d}(p_{k-1},
 p_k).$$

\definition{Definition 6.9} Let $\tau_{t_0}$ denote 
the time translation operator with translation 
amount $-t_0$, i.e. $\tau_{t_0}(u)(y, t)= u(y, t-t_0)$.
\enddefinition

\definition{Definition 6.10} We define 
the ``suspension" or ``pre-gluing" operator
$$
\#: L_{l, \d}(p_0,...p_k) \times \Bbb R^{k-1}_+
\to L_{l, \d}(p_0, p_k),
$$
$\Bbb R_+ = \{ t \in \Bbb R: t >0 \}$,  as follows.
Let 
$\chi$ be the cut-off function introduced 
before.
For $(u, \bold r) = ( (u_1, ..., u_k), (r_1, ...,
r_{k-1}))  \in L(p_0, ..., p_k) \times 
\Bbb R_+^{k-1}$, we set
$$
\#(u, \bold r) = 
\cases 
u_1, \hbox{ if } t \leq r_1-1; \cr
u_k, \hbox{ if } t \geq 2 \sum_{ i \leq k-2}
r_i + r_{k-1}+1; \cr
p_1 + \tau_{r_1}(\chi) (u_1-p_1) 
+ \tau_{r_1}(1-\chi) \tau_{2 r_1} (u_2-p_1), 
\hbox{ if } r_1-1  \leq t \leq r_1 \cr 
 \hskip 160pt +1; \cr
... \cr
p_{k-1} + \tau_{2r_1+...2r_{k-2}
+r_{k-1}}( 
\chi ) \tau_{2 r_1 +...2r_{k-2}                 
} (u_{k-1}-p_{k-1}) + \cr
\hskip 10pt \tau_{ 2r_1+... r_{k-2} +
r_{k-1}}(1-\chi) 
\tau_{2 r_1+...2r_{k-1}}(u_{k}-p_{k-1} ), 
\hbox{ if }  2\sum_{j \leq k-2} r_j
 + \cr \hskip 50pt r_{k-1}  -1  \leq t \leq 2 
  \sum_{j \leq k-2} r_j 
 + r_{k-1}+1. \endcases
 $$
 We also denote $\#(u, \bold r) $ by 
 $u \# \bold r$.
\enddefinition

\definition{Definition 6.10} For $\bold r=(r_1,...r_{k-1})
 \in \Bbb R^{k-1}_+$
we introduce the weight
$w_{\bold r} \equiv exp(w^*_{\bold r})$,
$$
w^*_{\bold r} = \tau_{r_1}(\chi) \d_F + 
\tau_{r_1}(1-\chi) \tau_{2r_1+r_2}(\chi) 
\tau_{2r_1}(\d_F)+...+
$$
$$
\tau_{2r_1+...2r_{k-2}+r_{k-1}}(1-\chi)
\tau_{2r_1+...2r_{k-1}}(\d_F).
$$
The weighted Sobolev norms $\| \cdot \|_{l, w_{\bold r}}$
are obtained by replacing $\d_F$ in $\| \cdot \|_{l, \d}$
with $w_{\bold r}$. For $\bar u =(u, \bold r)$, we 
set 
$$
w_{\bar u} = w_{\bold r}.
$$
\enddefinition

\definition{Definition 6.11} For 
$u \in L_{l,\d}(p, q) \cup L^{l+1, 2}_{loc}$
with $
p, q \in \A_l(Y) \times \Ga_l(Y)$
 and 
$\bold r = (r_1,...,r_{k-1})$ 
$\in \Bbb R_+^{k-1}$ we define 
the $\bold r$-interpolation of 
$u$ as follows
$$
Int_{\bold r}(u)= u_0 + 
\chi (p-u_0)+ 
\tau_{r_1}(\eta) ( u(\cdot, r_1) -u_0)
+...+ $$
$$\tau_{2r_1+...2r_{k-2}+r_{k-1}}(\eta) ( 
u(\cdot, 2r_1+...2r_{k-2}+ r_{k-1})-u_0)+
$$
$$
\tau_{2r_1+...2r_{k-1}}(1-\chi)(q-u_0),
$$
where $u_0$ is the reference element introduced 
in Section 4 and 
$\eta$ is a cut-off function satisfying 
$$
\eta(t)= \cases 1, & \hbox{ if } -1 \leq t \leq 1; \cr
0, & \hbox{ if } t \leq -2  \hbox{ or } t \geq 2.
\cr
\endcases
$$
\enddefinition

\definition{Definition 6.12} For $\bold r=
(r_1,...,r_{k-1}) \in \Bbb R_+^{k-1}$, 
$l \geq 0$ and $u \in L_{l, \d}(p_1, q_2) \cap
L^{l+1, 2}_{loc}, v \in L_{l, \d}(p_2, q_2) 
\cap L^{l+1, 2}_{loc}$
with $p_1, p_2, q_1, q_2 \in \A_l(Y) \times \Ga_l(Y)$, 
we introduce 
the $(l, \bold r)$-distance 
$$
d_{l, \bold r}(u, v)^2= \|p_1-
p_2\|_{l, 2}^2 + \|q_1
-q_2\|_{l, 2}^2 + 
$$
$$
\sum_{1\leq i \leq k-1}
 \|u(\cdot, r_{i})-u(\cdot, r_i)\|_{l, 2}^2 +
\|(u-Int_{\bold r}(u))-
(v- Int_{\bold r}(v)) \|_{l, w_{\bold r}}^2.
$$
For gauge equivalence classes $\omega_1, \omega_2$
we set
$$
d_{l, \bold r}(\omega_1, \omega_2)=
inf\{d_{l, \bold r}(u, v): u \in \omega_1, v \in \omega_2\}.
$$
\enddefinition

\definition{Definition  6.13} We say that a sequence 
$u_j \in  L_{l, \d}(\A_l(Y) \times \Ga_l(S), 
\A_l(Y) \times \Ga_l(S)) \cap L^{l+1, 2}_{loc}$ converges to $u=
(u^1, ...,u^k) \in L_{l, \d}(p_0, ..., p_k)$ 
in ``piecewise exponential Sobolev $(l, 2)$-norm", 
or  ``piecewise exponentially in Sobolev
$(l, 2)$-norm", provided 
that there is  a sequence $\bold r_j=(r_{j, 1},...,
r_{j, k-1} ) \in \Bbb R_+^{k-1}$ such that 
 $r_{j, i} \to \infty$ for each $i$ and 
$d_{l, \bold  r_j}(u_j, \#(u, \bold r_j))
\to 0$ as $j \to \infty$.

The essense of this convergence concept is
this: roughly speaking, we can divide each $u_j$ into 
$k$ portions  to produce $k$ new sequences,
such that each of them
converges in exponential Sobolev $(l, 2)$-norm. One important 
feature is the ``equal distance property" of 
the exponential weight $w_{\bold r}$. Namely 
it consists of $k-1$ weights, such that 
the distance from the $i$-th ``lowest position" 
$r_i$ to the ``center position" of the 
$i$-th weight is equal to its distance to 
the center position of the $(i+1)$-th weight.
(The picture is clear when one draws the graph of 
the weight $w_{\bold r}$.) 

Smooth piecewise exponential convergence is 
defined to be piecewise exponential convergence 
with respect to all Sobolev norms.

Note that the one piece case $k=1$ of piecewise
exponential convergence is identical to 
the concept of exponential convergence given 
in Definition 4.6.
\enddefinition

A consequence of the piecewise 
exponential convergence  is the 
energy convergence, namely
$$
E(u_j) \to E(u) \equiv \sum_i E(u^i).
$$

\definition{Definition 6.14} 
The convergence concepts in Defintion 6.13 are 
extended in an obvious
way to  convergence of 
elements in $L_{l, \d}(p_0, ..., p_k)$ (with 
a fixed $k$ but varying $p_i$ in general).
 Namely the convergence 
is defined to be that of each piece.

 The piecewise exponential 
 convergence of equivalence classes 
 $[u^j]$ (in 
various contexts and set-ups) 
 is defined in terms of  
the piecewise exponential convergence  of $\ti u^j$ for
 some  
representives $\ti u^j \in [u^j]$. Equivalently,
we can use $d_{l, \bold r}$ for equivalence classes.

Then the  concept of piecewise exponential 
convergence is  also  naturally extended to 
multiple temporal trajectories classes. 
It is defined to be the piecewise 
exponential convergence of each piece.
\enddefinition

The real line $\Bbb R$ acts on trajectories in terms of 
 the time translation.  We define the time translation 
 action on multiple  
 temporal trajectories classes 
 or piecewise trajectories to be the separate time translation 
 action on individual pieces of the  classes. It 
 gives rise to an action of $\Bbb R^k$. Let the
  underline denote quotient
under the time translation action in the 
{\it 
context of temporal trajectories}, e.g. 
${\underline \M}_T^0(S_{\a},
S_{\b})= \M_T^0(S_{\a}, S_{\b}) \slash\Bbb R$.

\proclaim{Theorem  6.15}
Assme that $(\pi, \lambda)$ is $Y$-generic. 
Then the moduli spaces $ {\underline  \MM}_T^0(S_{\a},$  
$ S_{\b})$ are compact Hausdorff spaces, where 
the topology is given by  
smooth piecewise exponential convergence. 
(This topology is equivalent to that given by  
piecewise exponential convergence in the 
Sobolev $(2, 2)$-norm.)  
Similarly, for $p, q \in \SS \W_0$, the quotient of $
{\MM}^0(p, q)$ under the time translation action are 
compact Hausdorff spaces.  
\endproclaim

The proof of this theorem will be presented in 
Part II.

Consider $p \in \a, q \in \b$
such that $p, q \in \SS\W_0$. The 
spaces $ \MM^0(p, q)$ (or $\M^0(p,q)$) and ${ 
\MM}_T^0(S_{\a}, S_{\b})$ (or $\M_T^0(S_{\a}, S_{\b})$) are 
isomorphic. But their quotients under the time translation 
action are not isomorphic. For our purpose, the time translation action 
on the temporal model is more suitable. 
For this reason, we consider the action of $\Bbb R$ on $\N(p, q)$ 
induced from the time translation action on $\N_T(S_{\a}, S_{\b})$, 
which we call the ``twisted time translation".

\definition{Definition 6.16}
The twisted time translation $ T_R$ by $R$ 
is defined as follows. Let $u \in \N(p, q)$. Then $ T_R u  
= (g_T(u)^{-1})^*(\tau_R(T_G(u)))$. (See the proof of Lemma 6.1 for 
$g_T$ and $T_G$.) The twisted time translation acts on each 
piece
of $k$-trajectories separately, giving rise to a $\Bbb R^k$ action.
We use the underline to denote the quotient under the twisted 
time translation action.  
\enddefinition

\proclaim{Lemma 6.17} Assume 
that $(\pi, \lambda)$ is $Y$-generic. Then
the moduli spaces  
$\underline {\M}^0(p_0,...,$
$p_k)$ and
 $\underline {\M}_T^0(S_{\a_0},
...,S_{\a_k})$ are canonically diffeomorphic 
smooth manifolds, where $p_i, \a_i$ are the same as in Lemma 6.8. 
\endproclaim

\proclaim{Proposition 6.18} Assume that 
$(\pi, \lambda)$ is $Y$-generic. Then 
the moduli 
spaces $\underline {\MM}^0(p, q;\SS \W_0)$ are compact
Hausdorff spaces,
where $p , q \in \SS\W_0, p \not = q$.
\endproclaim 

The proof of this proposition will be presented 
in Part II.

Notice that the compact space ${\underline {\MM}}_T^0(
S_{\a}, S_{\b})$ contains ${\underline \M}_T^0(S_{\a},
S_{\b})$ as a subspace. The compactification of
 ${\underline \M}_T^0(S_{\a},
S_{\b})$ is given by its closure 
${\underline {\bar \M}}_T^0(S_{\a}, S_{\b})$. 

\proclaim{Theorem 6.19} For generic $(\pi, \lambda)$,
we have $\underline {\bar \M}_T^0(S_{\a},
S_{\b}) = {\underline {\MM}}_T^0(S_{\a}, S_{\b})$.
Moreover, the following hold:

(1) ${\underline {\MM}}_T^0(S_{\a},S_{\b})$ has the
 structure of 
$d$-dimensional smooth orientable  manifolds with 
corners (i.e. 
modelled on the first quadrant of $\Bbb R^d$), where 
$d = \mu(\a)-\mu(\b)-\hbox{dim }G_p +1$.

(2) This structure is compatible with the stratification 
${\underline {\MM}}_T^0(S_{\a}, 
S_{\b})= \cup_k {\underline { \M}}_T^0(S_{\a},$
$S_{\b})_k$, i.e. the interior of the $k$-dimensional face of 
${\underline {\MM}}_T^0(S_{\a}, S_{\b})$ 
is precisely ${\underline {\MM}}_T^0(S_{\a},$ $S_{\b})_k$.

(3) The temporal projections $\pi_-:
{\underline {\MM}}_T^0(S_{\a},
S_{\b}) \to  S_{\a}$ and $\pi_+: {\underline {\MM}}_T^0(S_{\a},
S_{\b})$ $ \to S_{\b}$ are smooth fibrations.

Consequently, we have 
$$\eqalign{
\pa \underline {\MM}_T^0(S_{\a}, S_{\b})&=
\cup_{\mu(\a) > \mu(\gamma) > \mu(\b)}
\underline {\MM}_T^0(S_{\a}, S_{\ga},
S_{\b})\cr
&\equiv \cup_{\mu(\a) > \mu(\gamma) > \mu(\b)}
\underline {\MM}_T^0(S_{\a}, S_{\ga}) 
\times_{S_{\ga}} \underline {\MM}_T^0(
S_{\ga}, S_{\b}),\cr} \tag 6.1 
$$
where the fiber product space $\underline {\MM}_T^0({\a}, S_{\ga}) \times_{S_{\ga}}
\underline {\MM}_T^0(S_{\ga}, S_{\b})$
is defined to be 
$$\{(u, v)\in \underline {\MM}_T^0(S_{\a}, S_{\ga}) \times
\underline {\MM}_T^0(S_{\ga}, S_{\b}):
\pi_+(u) = \pi_-(v)\}. 
$$
(Of course, only the nonempty moduli spaces 
appear in (6.1).)
\endproclaim

The situation with the isomorphic model 
$\M^0(p, q; \SS\W_0)$ is similar. Namely, the compact 
space ${\underline {\MM}}^0(p, q; \SS\W_0)$
contains the subspace $\M^0(p, q; \SS\W_0)$, whose 
closure we denote by ${\underline {\bar \M}}^0(p, q;\SS\W_0)$.
Then we have the following equivalent formulation for 
Theorem 6.19. (We only give part of the statements.)

\proclaim{Theorem 6.20} In a generic situation, we have 
${\underline {\bar \M}}^0(p, q; 
\SS\W_0)= \underline {\MM}^0(p, q; $
$\SS\W_0)$.
Moreover, ${\underline {\MM}}^0(p,q;\SS\W_0)$ has the 
structure of 
$d$-dimensional manifolds with corners which is 
compatible with its natural stratification.
\endproclaim 

The said equivalence means the following:
\proclaim{Proposition 6.21}
In a generic situation, ${\underline {\MM}}^0(p, q)$
is canonically diffeomorphic to ${\underline
 {\MM}}_T^0( S_{\a}, S_{\b})$.  
\endproclaim

This proposition is an easy consequence of the temporal 
transformation. 
Clearly, Theorem 6.20 is a consequence of 
Theorem 6.19.
The proof of Theorem 6.19  will be presented in
Part II.

The two different models- the 
{\it temorpal model } $\underline \MM^0_T(S_{\a},
S_{\b})$ and the {\it 
fixed-end model} $\underline \MM^0(p, q; \SS\W_0)$- provide 
different aspects of our basic set-up.  It is good to have 
both of them for the purpose of conceptual understanding.
Since the former is canonical (while the latter involves
the choice of $\SS\W_0$), we shall use it in the 
formulation of our 
main constructions in the sequel. Of course, we can 
use equally well the latter.

\head 7. Bott-type homology
\endhead

We first introduce a few orientation conventions. We 
follow those used in 
$\FukayaI$.  For an oriented smooth manifold with 
corners $\X$, its boundary is oriented in 
such a way that $span\{n_{\pa \X}\} \oplus
T\pa \X = TX|_{\pa \X}$ as oriented vector bundles (away from the corners 
of $\pa \X$), where $n_{\pa \X}$ is an inward normal field of the 
boundary. Given transversal smooth maps 
$F_1: \X_1 \to S$ and $F_2: \X_2 \to S $ from 
two oriented smooth manifolds with corners into 
an oriented smooth manifold $S$, the  fibered 
product $\X_1 \times_S \X_2 = (F_1 \times F_2)^{-1}
(Diag (S \times S))$ ($Diag$ means {\it diagonal})
has a canonical orientation such that $T(\X_1 \times_S 
\X_2) \oplus N = (-1)^{dim S \cdot dim \X_2} TX_1 \oplus TX_2$
as oriented bundles, where $N$ denotes the 
oriented bundle $(d (F_1 \times F_2))^{-1}((TS \oplus \{0\})|_{Diag(
S \times S)})$. The following lemma can 
be found in $\FukayaI$  and is easy to verify. 

\proclaim{Lemma 7.1} 
There holds for oriented boundaries
$$
\pa(\X_1 \times_S \X_2)= \pa \X_1
\times_S \X_2 + (-1)^{ dim  \X_1
+  dim  S} \X_1 \times_S \pa \X_2, \tag 7.1
$$
where the summation sign means taking union.
We also have
$$ (X_1 \times_{S} \X_2) \times_{S'} \X_3 =
X_1 \times_S (\X_2 \times_{S'} \X_3). \tag 7.2
$$
\endproclaim 

We shall use the natural orientation of $S_{\a}$ 
induced from the $S^1$ action. 
(We can also use any other orientations.)
By Theorem 6.13, the boundary of 
$\underline {\MM}_T^0(S_{\a}, S_{\b})$ is 
a union of fibered products. We need to arrange the orientation 
of these spaces so that a suitable consistency holds in regard 
of the natural orientation of fibered products as defined 
above and boundary orientations. Indeed, we have 

\proclaim{Lemma 7.2} We can choose the orientation of 
$\underline {\MM}_T^0(S_{\a}, S_{\b})$ such that
$$
\pa \underline {\MM}_T^0(S_{\a}, S_{\b}) = 
(-1)^{\mu(\a) + dim S_{\a}} \Sigma_{\mu(\a)  > \mu(\ga) 
> \mu(\b)} \underline {\MM}_T^0(S_{\a}, 
S_{\ga}) \times_{S_{\ga}} \underline {\MM}_T^0(S_{\ga},
S_{\b})
\tag 7.3
$$
as oriented manifolds with corners.
\endproclaim  

This lemma is analogous to Sublemma 1.20 in $\FukayaI$ and can be  proven 
by the same arguments as there.

\bigskip


\bigskip

For a topological space $\X$ and an abelian group
$\GG$, let $ Q_j(\X; \GG)=Q_j(\X; \Bbb Z) \otimes \GG$
 denote
the abelian group   generated by generalized singular
$j$-cubes in $\X$ with coefficients in  $\GG$.
Here we define a generalized  
 singular $j$-cube to be  
a continuous map $f: \Delta \to \X$, where 
$\Delta$ is an oriented  $j$-cube, namely 
an oriented smooth manifold with corners which is 
diffeomorphic to the euclidean $j$-cube
$[0, 1]^j$.
(A continuous map $f: [0, 1]^j \to 
\X$ is called a singular $j$-cube, see 
$\Massey$.)
 Let $D_j(\X;
\GG)$ be the  subgroup  of $Q_j(\X; \GG)$
generated by degenerate generalized singular $j$-cubes. Here, 
a generalized singular $j$-cube $f: \Delta \to \X$  is called 
degenerate, if $f \circ F$  is independent of some variable
$t_i$, where $F$ is a diffeomorphism from 
$[0, 1]^j$ to $\Delta$.   

Next consider a 
$j$-dimensional compact oriented manifold
with corners $\Sigma$ and a continous map
$f: \Sigma \to \X$. We divide  $\Sigma$ into
oriented cubes and produce this way  
a chain $\sigma \in Q_j(\X; \Bbb Z)$.
Choosing a different  way of dividing,
we obtain another chain $\sigma'$. 
Let $Div_j(\X; \GG)$ be the 
subgroup of $Q_j(\X; \GG)$ generated by 
elements of  the form $(\sigma-\sigma') 
\otimes g , g \in \GG$.  

 Let $ C_j(\X; \GG)$ denote 
 the quotient group $Q_j(\X; \GG)/(D_j(\X; \GG)
 +Div_j(\X; \GG))$. This is  
 the group of generalized cubical singular 
 $j$-chains in $\X$ (with coefficients in
 $\GG$).  The corresponding 
 cochain group is $C^j(\X; \GG)=Hom( C_j(\X;
 \Bbb Z), \GG)$.

 The ordinary boundary operator $\partial_{O}$ 
 for singular $j$-cubes extend to generalized 
 singular $j$-cubes canonically via 
 orientation preserving diffeomorphisms 
 between $j$-cubes and $[0,1]^j$.  It then extends 
 to $Q_j(\X; \GG)$ by linearity, and to $ C_j(\X;
 \GG)$ by descending.    We 
have the corresponding coboundary operator $\pa_O^*$.
 It is easy to show that the homology $H_*( 
 C_*(\X; \GG), \pa_O)$ and cohomology $H^*( 
 C^*(\X; \GG), \pa_O^*)$ are canonially isomorphic to 
 the singular homology and cohomlogy of $\X$ 
 with coefficient group $\GG$.

\definition{Remark 7.3} For $\sigma= ([0,1]^j, f) \in C_j(\X;
\Bbb Z)$,
the map $( \pa [0, 1]^j, f|_{\pa [0, 1]^j})$
from the oriented boundary induces a chain in 
a natural way. By our convention for boundary 
orientation, this chain equals $- \pa_O \sigma$.
\enddefinition

Now we fix a coefficient group $\GG$ and abbreviate 
e.g. $ C_j(\X; \GG)$ to $ C_j(\X)$.
We set 
$$S_i=\cup \{S_{\a}: \a \in\R, \mu(\a)=i\}
\tag 7.4
$$
and 
introduce our Bott complexes $C_*^{Bott}$ and
$C^*_{Bott}$,
 $$C_k^{Bott}=\oplus_{i+j=k}C_j(S_i),
C^k_{Bott}=\oplus_{i+j=k}C^j(S_i).$$
Note 
$ C_*^{Bott} =  C_*(\R^0), C^*_{Bott}
=   C^*(\R^0). $
                      
Next we define a boundary operator
$\pa_{Bott} : C_k^{Bott} \to C_{k-1}^{Bott}$ 
along with a coboundary operator 
$\pa^*_{Bott}: C^{k-1}_{Bott} \to C^k_{Bott}$ for each $k$. 
First, for each pair $\a, \b \in \R$ with $\mu(\a)
> \mu(\b)$, 
 we define a boundary operator 
$\pa_{\a, \b}: C_k^{Bott} \to C_{k-1}^{Bott}$.
If the moduli space $\underline 
{\MM}_T^0(S_{\a}, S_{\b})$ is empty, we define 
$\pa_{\a, \b}$  
to be the zero operator. If it is nonempty, we define 
$\pa_{\a, \b}$ as follows. For $\sigma 
\not \in C_*(S_{\a})$, we set $\pa_{\a, 
\b} \sigma = 0$. For $\sigma=[(\Delta,f)] \in C_j(S_{\a})$
with $\mu(\a)+j=k$, 
consider the fibered product 
$$\ti \Delta=
 \Delta \times_{S_{\a}}\underline {\MM}_T^0(S_{\a},
 S_{\b}) =\{(z,u)\in \Delta\times
 \underline {\MM}_T^0(S_{\a},S_{\b}):f(z)=\pi_-(u)\}.$$

We have a natural 
projection map $ \pi_+: \ti \Delta \to 
S_{\b}, \pi_+((z, u))=\pi_+(u).$
Since the projection $\pi_-$ is a submersion, $\ti \Delta$ is 
a compact oriented manifold with corners of 
dimension $j+\mu(\a)-\mu(\b)-1$. 
We divide $\ti \Delta$ 
 into oriented cubes. Using the projection map
 $\pi_+$ we then obtain from the divided 
 $\ti \Delta$  an element    $\ti \Delta_f$
 in $Q_{j+\mu(\a)-\mu(\b)-1}(S_{\b}) $. We 
 define 
 $$
 \pa_{\a, \b}(\sigma) = [\ti \Delta_f].
 $$
 It is easy to see that this definition is 
 independent of the choice of $(\Delta, f)$
 and the way of dividing $\ti \Delta$.

 Clearly, we indeed have
 $$
 \pa_{\a, \b}: C_k^{Bott} \to C_{k-1}^{Bott}
 $$
for all $k$.

\definition{Definition 7.4} We define 
$ \pa_{Bott} : C_k^{bott}
\to C_{k-1}^{Bott}$ as follows:
$$
 \pa_{Bott} = \pa_0  + 
\sum_{\mu(\a) > \mu(\b)} \pa_{\a, \b},$$ 
where  
$ \pa_0= (-1)^k \pa_{O}$. 
The 
boundary operator $\pa_{Bott}: C_*^{Bott} \to
C_*^{Bott}$ is defined to be the direct sum of these
 boundary operators.
  
The boundary operator $\pa_{Bott}$ induces a coboundary 
operator $\pa^*_{Bott}: C^k_{Bott} \to C^{k+1}_{Bott}$.
Indeed, we have
 $\pa^*_{Bott} =(-1)^k \pa_O^* + \sum \pa^*_{\a, \b}$.
\enddefinition

\proclaim{Lemma 7.5} We have $\pa_{Bott}^2=0, (\pa^*_{Bott})^2
=0.$ 
Hence $(C_*^{Bott}, \pa_{Bott})$ is a chain 
complex and $(C^*_{Bott}, 
\pa^*_{Bott})$ is a cochain complex.
 \endproclaim

 \demo{Proof}
We handle the case $\GG=\Bbb Z$ and $\pa_{Bott}$. The arguments 
easily extend to $\pa^*_{Bott}$ and general coefficient 
groups.
Consider $\sigma = [(\Delta, f)]\in C_j(S_{\a})$
with $j+ \mu(\a)=k$. 
 We have $\pa_{Bott}^2 \sigma = \sum_{\mu(\a) <\mu(\b)}
  I_{\b}$, where 
 $$
 I_{\b} = \pa_0 \pa_{\a, \b}\sigma
 + \pa_{\a, \b}\pa_0\sigma
 + \sum_{\mu(\a)> \mu(\gamma) > 
 \mu(\b)}\pa_{\gamma, \b} \pa_{\a,
 \gamma} \sigma.  
$$
On the other hand, by (7.1) we have
$$\pa(\D  \times_{S_{\a}} \underline 
{\MM}_T^0(S_{\a},S_{\b})) =
\pa \D \times_{S_{\a}}\underline {\MM}_T^0(S_{\a},S_{\b})
+ 
(-1)^{j+dim S_{\a}}\D_j \times_{S_{\a}}\pa \underline 
{\MM}_T^0(S_{\a},S_{\b}) \tag 7.5 $$ 
$$= 
\pa \D \times_{S_{\a}}\underline {\MM}_{T}^0(S_{\a},S_{\b})
+ 
$$
$$(-1)^{j+ \mu({\a})}\sum_{\mu(\a) 
>\mu(\ga) > \mu(\b)} \D \times_{S_{\a}} \underline 
{\MM}_T^0(S_{\a},S_{\ga}) \times_{S_{\ga}}
\underline {\MM}_T^0(S_{\ga},S_{\b}).
$$
Multiplying this equation by
 $(-1)^{j+ \mu(\a) -1}=(-1)^{k-1}$ 
 and applying Definition 7.4 we then deduce
 $$
 \pa_0 \pa_{\a, \b} \sigma =- \pa_{\a, \b} \pa_0 \sigma
 -\sum_{\mu(\a) > \mu(\gamma) > \mu(\b)} \pa_{\a,
 \gamma} \pa_{\gamma, \b} \sigma.
 \tag 7.6
 $$
 This shows that $I_{\b} =0$. 
\qed 
\enddemo

\definition{Definition 7.6} We define the Bott-type Seiberg-Witten 
Floer homology $$FH^{SW}_{Bott; \GG}(c)_*$$
to be the 
homology $H_*(C_*^{Bott}, \pa_{Bott})$ of the 
chain complex $(C_*^{Bott}, \pa_{Bott})$.
 (Recall that 
$c$ is the given $spin^c$ structure.) We define 
the Bott-type Seiberg-Witten Floer cohomology 
$FH^{SW}_{Bott; \GG}(c)^*$ to be the cohomology 
$H^*(C^*_{Bott}, \pa^*_{Bott})$ of the 
cochain complex $(C^*_{Bott}, \pa^*_{Bott})$.
\enddefinition

Next we introduce a natural filtration 
-the index filtration $$\F_*^{Bott}= \cdot \cdot \cdot \subset
\F_k^{Bott} \subset \F_{k+1}^{Bott} 
\subset \cdot \cdot \cdot $$ 
on our chain complex $(C^{Bott}_*, \pa_{Bott})$.
We put
$$
\F_k^{Bott} = \oplus_{m \leq k}  C_*(S_m).
$$
It is clear that this defines a filtration by subcomplexes.
We also have a dual filtration $\F^*_{Bott} =
\cdot \cdot \cdot \F^k_{Bott} \subset \F^{k-1}_{Bott}
 \subset
\cdot \cdot \cdot$ for the cochain complex,
$$
\F^k_{Bott}= \oplus_{m \geq k} C^*_m.
$$

\proclaim{Theorem 7.7} The index filtration induces
 a spectral sequence
$E(Bott)^*_{**}$ which converges to
$FH^{SW}_{Bott*}$
and satisfies 
$$
E(Bott)^1_{ij} \cong H_j(S_i, \GG). \tag 7.7
$$
The dual filtration induces a  spectral sequence 
$E(Bott, dual)^*_{**}$ which converges to $FH^{SW*}_{Bott}$
and satisfies 
$$   E(Bott, dual)^1_{ij}
\cong  H^j(S_i; \GG). \tag 7.8
$$
\endproclaim

\demo{Proof} The construction of the spectral sequence out of 
the filtration is 
standard. The  associated graded complex $G C_*^{Bott}$ for our 
filtered complex $(C_*^{Bott}, \F_*^{Bott}) $ is given by
$$
GC_*^{Bott} = \oplus_k C_*({S_k}),
$$
which leads to the formula (7.8). 
The convergence of the spectral sequence to the Bott-type 
Seiberg-Witten Floer homology  is a consequence of the 
basic theory of spectral sequences.   The case of the dual filtration 
is similar.
\qed
\enddemo

Finally, we introduce a tool which will 
be used for proving the invariance of 
our Seiberg-Witten Floer homology theories.

\definition{Definition 7.8}
Consider a finite collection  $\bold F=\{f_1,...,f_m\}$, 
where each $f_i$ is a smooth map from a compact 
smooth manifold with corners to some $S_{\a}$.
A generalized singular cube $f: \Delta \to S_{\a}$ for 
some $\a$ is called $\bold F$-transversal, provided 
that the following holds: 

\noindent 1) $f$ is smooth,

\noindent 2) $f$ and $f|_{
\partial \Delta}$ are transversal  to every 
$f_i$,  and

\noindent 3) for each 
sequence $\b_1, ..., \b_j$ with 
$\mu(\a) > \mu(\b_1) >...> \mu(\b_j)$, 
the induced maps $\pi_+: 
\Delta \times_{S_{\a}} (\MM_T^0(S_{\a},
S_{\b_1}) 
\times_{S_{\b_1}} ... \times_{S_{\b_{j-1}}}
\MM_T^0(S_{\b_{j-1}}, S_{\b_j})) \to 
S_{\b_j}$  and $\pi_+:
\pa \Delta \times_{S_{\a}} (\MM_T^0(S_{\a},
S_{\b_1}) 
\times_{S_{\b_1}} ... \times_{S_{\b_{j-1}}}
\MM_T^0(S_{\b_{j-1}}, S_{\b_j}) )\to 
S_{\b_j}$ are transversal to every 
$f_i$.

We have
the subgroup $Q_j^{\bold F}(\R^0)$ 
of $Q_j(\R^0)$ generated by 
$\bold F$-transversal generalized singular cubes,
which yields a subgroup $C_j^{\bold F}(\R^0)$
of $C_j(\R^0)$, whose elements will be called
$\bold F$-transversal chains.
We define 
$$
C_k^{Bott, \bold F}= \oplus_{i+j=k} C_j^{\bold F}
(S_i)$$
and obtain therewith
a subcomplex $C_*^{Bott, \bold F}$
of $C_*^{Bott}$. 
It is easy to see that this is indeed a 
subcomplex. We have the corresponding cochain complex 
$C^*_{Bott, \bold F}$ dual to $C_*^{Bott, 
\bold F}$.
\enddefinition

For each given $\bold F$, we  obtain a homology
$H_*(C_*^{Bott, \bold F}, \pa_{Bott})$
and a  cohomology $H^*(C^*_{Bott, \bold F},
\pa_{Bott}^*)$.

\proclaim{Lemma 7.8} The homology 
$H_*(C_*^{Bott, \bold F}, \pa_{Bott})$ is 
canonically isomorphic to the homology 
$H_*(C_*^{Bott},
\pa_{Bott})$. Similarly, the cohomology 
$H^*(C^*_{Bott, \bold F}, \pa^*_{Bott})$ 
is canonically isomorphic to the cohomology
$H^*(C^*_{Bott}, 
\pa_{Bott}^*)$.
\endproclaim
\demo{Proof} We have the inclusion 
homomorphism 
$I: C_*^{Bott, \bold F} \to 
C_*^{Bott}$ and the induced homomorphism $I_*:
H_*(C_*^{Bott, \bold F}, 
\pa_{Bott}) \to H_*(C_*^{Bott}, \pa_{Bott})$.
We claim that $I_*$ is an isomorphism. 
First, consider a closed chain $\sigma $ 
in $C_*^{Bott, \bold F}$ such that 
$F(\sigma)= \pa_{Bott} \sigma_1$ 
with $\sigma_1 \in C_*^{Bott}$. 
By smooth approximation 
and transversal perturbations we 
deform $\sigma_1$ to obtain $\sigma_2 
\in C_*^{Bott, \bold F}$ such that 
$\sigma = \pa_{Bott} \sigma_2$. 
This shows that $F_*$ is a monomorphism.
On the other hand, we can use smooth 
approximation and transversal perturbations 
to deform an arbitary 
closed chain $\sigma \in C_*^{Bott}$ 
to obtain a closed chain $\sigma' 
\in C_*^{Bott, \bold F}$ 
such that $F(\sigma')- \sigma$ is 
homologeous to zero. This shows that 
$F_*$ is an epimorphism. 

The arguments for the cohomologies are 
similar.
\qed
\enddemo

\head  8. Invariance 
\endhead

Consider two chain complexes $(C_*, \pa)$ and 
$(\bar C_*, \bar \pa)$. A chain map  of degree $m \in \Bbb 
Z$  from the former to the latter  
consists of homomorphisms $F: C_k \to \bar C_{k+m}$ 
such that $F \cdot \pa = \bar \pa \cdot F$. 
A shifting homomorphism $F: H_*(C_*, \pa) 
\to H_*(\bar C_*, \bar \pa)$ of degree $m \in \Bbb Z$ consists 
of homomorphisms $F: H_k(C_*,\pa) \to H_{k+m}(\bar C_*, \bar \pa)$. 
These concepts are also defined 
for cochain complexes and cohomologies  
in a similar way, with the degree shifting in the 
opposite way, i.e.
$k \to k-m$.

We shall establish the following invariance result.

\proclaim{ Theorem 8.1} 
The Bott-type Seiberg-Witten Floer 
homology and cohomology are diffeomorphism 
invariants  modulo shifting isomorphisms. 
\endproclaim 

In othere words, these homology and cohomology are independent of 
the metric $h$ and generic parameter $(\pi, \lmd)$
modulo shifting isomorphisms.

In this section, we construct a shifting 
homomorphism, which will be shown in  $\YeII$ 
to be a desired shifting ismorphism. 
We would like to point out that the construction 
already contains some basic arguments for 
establishing the  
isomorphism property.

Consider two metrics $h_+$ and $h_-$ on $Y$ and 
$Y$-generic parameters $(\pi_+,\lmd_+)$ for $h_+$ and 
$(\pi_-, \lmd_-)$ for $h_-$ respectively. We would like to 
construct a shifting isomorphism between  our 
homologies (cohomologies) constructed with respect to $(h_+, \pi_+,
\lmd_+)$ and $(h_-, \pi_-, \lmd_-)$
respectively. 
 
Choose a smooth path of metrics $h(t)$ on  $Y$  such that
$$h(t)=\cases h_-, & \hbox{if }t<-1,\cr
h_+, & \hbox{if }t>1,\cr\endcases
$$
a smooth path of $\pi(t)\in\Pi$ 
$$\pi(t)=\cases \pi_-, & \hbox{if }t<-1,\cr
\pi_+, & \hbox{if }t>1,\cr\endcases
$$
and a smooth function $\lmd(t)\in \Bbb R$ such that
$$\lmd=\cases \lmd_-, & \hbox{if }t<-1,\cr
\lmd_+, & \hbox{if }t>1.\cr\endcases
$$

The following lemma is an immediate consequence of 
Lemma B.1 in Appendix B.

\proclaim{Lemma 8.2} Fix an $(A_0, \P_0)
\in \A(X) \times \Ga^+(X)$. Indeed,
we use the reference element $u_0$ introduced 
before as $(A_0, \P_0)$.  For any $R>0$, set $X_R= Y \times 
[-R, R]$ and
$$\eqalign{\SS_R=\{&u = (A_0, \P_0)+ (A, \P): (A, \P) 
\in \O^1_2
(X_R) \times \Ga^+_2(X_R)\cr
&\hbox{with }d^* A=0\hbox{ and }A|_{\pa X_R}(\frac{\pa}
{\pa t}) =0. \}\cr}$$
Then $\SS_R$ is a global slice for the action 
of $\G_3^0(X_R)$
on $\A_2(X_R) \times \Ga_2^+(X_R)$. In other words, $\A_2(X_R)
\times \Ga_2^+(X_R) = \G^0_3(X_R)\cdot \SS_R$.
\endproclaim

The following definition is a crucial construction. 

\definition{Definition 8.3} Choose a nonzero 
$\Psi_0 \in \Ga^-(X)$ with 
support contained in the interior of $X_1$. We define 
a smooth vector field $Z$ on $\A_{2}(X_1) \times 
\Ga^+_{2}(X_1) $ by
$$
Z(g^*u )=g^{-1}\Psi_0,
$$
for $g \in \G^0_3(X_1), u \in \SS_1$. Furthermore, 
we extend $Z$  to 
$\A_{2, loc}(X) 
\times \Ga^+_{2, loc}(X)$ by setting
$$
Z(u)=Z(u|_{X_1}).
$$
\enddefinition

The following lemma is readily proved. 
\proclaim{Lemma 8.4}
$Z$ is equivariant with respect to the action of $\G_{3,
 loc}^0$. 
\endproclaim

Let $X$ be endowed with 
the warped product metric 
determined by the family of metrics $h(t)$ and the standard
 metric 
on $\Bbb R$. We introduce the (perturbed) transition trajectory 
equation 
$$
\cases
F^+_A =&\frac{1}{4}\langle e_ie_j\P,\P\rangle e^i\wedge e^j
\cr 
&+ ( \nabla H_{\pi(t)}
(a) +b_0) \wedge dt+*(  (\nabla H_{\pi(t)} (a)+
b_0) \wedge dt),\cr
D_A\P =& -\frac{\pa}{\pa t} \cdot (\lambda \P- Z(A, \P)),\cr
\endcases
\tag 8.1
$$
where $b_0$ denotes a smooth 1-form of compact 
support on $X$ without $dt$-part.
Setting  $A= a + f dt, \P=\p$ as usual we can 
rewrite (8.1) as follows 
$$
\cases
\frac{\partial a}{\partial t} & =*_YF_a+d_Yf+\langle e_i\cdot\p,\p\rangle e^i+
\n H_{\pi(t)}(a) +  b_0,
\cr
\frac{\partial \p}{\partial t} & =-\np_a\p-\lmd(t)\p+ 
Z.\cr
\endcases
\tag 8.2
$$ 
Of course, in those places in (8.2) where 
 the metric on $Y$ is involved, we use the metric 
 $h(t)$. For example,
 the  Hodge $*_Y$ at time $t$ in the 
equation is that of the metric $h(t)$.  
The perturbation term $  Z$ (or $ \pa/ \pa t \cdot Z$) 
is called a {\it spinor 
perturbation}.  We have  the following obvious, but crucial lemma. 

\proclaim{Lemma 8.5} The equation (8.1) 
is invariant  with respect to the action of $\G^0_{3, loc}$.
Moreover, it has no reducible solution.
\endproclaim

A fundamental property of the Seiberg-Witten equation 
is a pointwise maximum principle  for the spinor 
field, which is a consequence of 
the Weitzenb{\"o}ck formula (2.1). With the presence of $Z$, this 
principle no longer holds. Instead, we have the following 
result.

\proclaim{Lemma 8.6} Let $(A, \P)=(a+ f dt, \p) $ be 
a solution of (8.1). Then there holds
$$
\|\P\|_{L^{\infty}} \leq C (1+E(A, \P))^{1/4}
$$
for a constant $C$ depending only on the 
families $h(t), \pi(t), \lmd(t)$ and the geometry of $Y$. 
\endproclaim
\demo{Proof} Before proceeding 
with the proof, we first observe that 
by (2.17) and (8.1) the energy can be estimated in the following way
$$
\eqalign {E((A,\P))=2 \lim_{t \to
-\infty} \hbox{ \bf cs}_{(\lambda, H)}(a(\cdot, t), \p(\cdot, t))
\cr -2\lim_{t \to +\infty} \hbox{\bf cs}_{(\lambda, H)}(\a(\cdot,
t), \p(\cdot, t))+\int_X|Z|^2, \cr} \tag 8.3$$ 
$$\int_X|Z|^2<C. \tag 8.4 $$
Using local Columb gauges provided 
by Lemma B.1 in Appendix B and a patching argument,
we can perform a gauge transformation to convert $(A, \P)$ into a smooth 
solution. Since the $L^{\infty}$ norm of $\P$ is invariant, 
we can assume that $(A, \P)$ is already smooth. Furthermore, 
we can assume that $(A, \P)$ is in temporal form.

For simplicity, we assume $\lmd(t) \equiv 0$ in the following 
argument. It is easy to modify the argument
 to handle $\lmd(t)$. Put
$$I_1= \np_a\p, I_2= *_Y F_a+\langle e_i\cdot\p,\p\rangle
e^i - \n  H(a, \p), I_3=I_2 + \n H(a, \p).$$
Then 
$$\int_{\Omega}(|I_1|^2+|I_2|^2) = \frac 12 E((A, \P), \Omega),$$
for any domain $\Omega \subset X$.
For each $t\in\Bbb R$, we use the 3-dimensional  Weitzenb\"ock formula (2.1) 
on $Y\times\{t\}$ to derive 
$$\np_aI_1=-\D_a\p+\frac{\bar s}4\p+\frac 12|\p|^2\p+
I_3\p,\tag 8.5$$
where $\bar s$ denotes the scalar curvature function of $(Y, h(t))$.
Multiplying (8.5) by $\p$ and integrating by parts, we infer  
$$\int_{Y\times\{t\}}(|\n_a\p|^2+\frac s4|\p|^2+ \frac 12 |\p|^4) \leq  
\int_{Y\times\{t\}}|I_1|\cdot|\pa_a\p|+|I_3||\p|^2.$$
Using the H\"older ineqality we then  deduce that
$$\int_{Y\times\{t\}}(|\n_a\p|^2+ \frac 14
|\p|^4) \leq C(1+\int_{Y\times\{t\}}(|I_1|^2+|I_2|^2),$$
where $C$ depends only on $\|\n H\|_{L^{\infty}}$, which 
can be estimated by appealing to Lemma 3.13.
This last estimate implies 
$$\int_{X_{R-2, R+2}}(|\n_a \p|^2+ |\p|^4) \leq C(1 + E(A, \P, 
X_{R-2, R+2}))
\tag 8.6$$
for any $R >0$, where $X_{r, R} =Y \times [r, R]$. 

Next we apply the Moser iteration to deduce 
the desired $L^{\infty}$ estimate. Let  
$\xi:X\to [0,1]$ be a cut-off function such that
${supp~}\xi \subset X_{R-2, R+2}$ and $\xi(t)=1$ 
for $t\in X_{R-1, R+1}$.  
By the 4-dimensional  Weitzenb\"ock formula (2.1) on $X$ (recall 
that $X$ is endowed with warped product metric), we have
$$
D_AZ=-\D_A \P+\frac{s}4 \P-\frac14 |\P|^2\P.
\tag 8.7$$
Choosing $\xi^2|\P|^p\P$ as a test function,
where $p>0$ will be determined later, we obtain 
$$\eqalign{\int_X\xi^2(\n_A\P\cdot\n_A(|\P|^p\P)+
\frac{s}4|\P|^{p+2}+\frac14 |\P|^{p+4})\cr
=-\int_X(\xi^2 ZD_A(|\P|^p\P)+2\xi\n \xi Z|\P|^p\P+
2\xi\n \xi\cdot\n \P|\P|^p\P).\cr}
\tag 8.8$$

We have 
$$
\eqalign{\n_A(|\P|^p\P) &=|\P|^p\n_A\P+d|\P|^p\P \cr
&=|\P|^p\n_A \P+p|\P|^{p-1}\frac{\langle \n_A\P,\P\rangle}{|\P|}\P,
\n_A\P\cdot\n_A(|\P|^p\P) \cr &=|\P|^p|\n_A\P|^2+p|\P|^{p-2}\langle\n_A\P,
\P\rangle^2. \cr}$$
On  the other hand,
$$ p |\P|^{p-2}  \langle \n_A \P, \P 
 \rangle^2 \leq C \varepsilon
|\P|^p |\n_A \P|^2 + \frac {Cp^2}{\varepsilon}  |\P|^p,
$$
where $\varepsilon > 0$ is arbitrary. Similarly,
$$\eqalign{|D_A(|\P|^p\P)| & \le|\P|^p|D_A\P|+|d|\P|^p\cdot\P| \cr
& \le C |\P|^p|\n_A\P|+p|\P|^p|\n_A\P| \cr
& \le C \varepsilon  |\P|^p|\n_A\P|^2+\frac C{\varepsilon} (p^2+1) |\P|^p,
\cr}$$

$$|\n \xi||\xi||Z||\P|^{p+1}\le C\xi^2|\P|^{p}
+C|\n\xi|^2|\P|^{p+2},$$
and
$$|\n \xi||\xi||\n
\P||\P|^{p+1}\le\varepsilon\xi^2|\P|^p |\n\P|^2+ \frac C{\varepsilon}
|\n \xi|^2|\P|^{p+2}.$$
Choosing $\varepsilon$ suitably, we deduce 
$$\int\xi^2|\P|^p|\n_A\P|^2\le C\int((p^2 +1)\xi^2|\P|^p+|\n
\xi|^2|\P|^{p+2}).$$
Consequently,
$$\int_X \xi^2|\P|^p|\n|\P||^2\le C(\int_X(p^2+1) \xi^2|\P|^p+
|\n \xi|^2|\P|^{p+2}),$$
or
$$
\int_X\xi^2|\n|\P|^{\frac{p+2}2}|^2\le C(p+1)^2\int_X(\xi^2|\P|^p+|\n
\xi|^2|\P|^{p+2}).
$$
Now we set $w=|\P|^{(p+2)\slash 2}$.
By the Sobolev inequality and  H\"older inequality, we arrive at 
$$\eqalign{\|\xi w\|_{L^2}^2 & \le C(p+1)^2\int 
(|\xi\n w|^2+|w\n \xi|^2)\cr
&\le C(p+1)^2(\|\xi
w\|_{L^2}^{\frac{p}{p+2}}+\||\n\xi|w\|_{L^2}).\cr}$$
Then we use the iteration process as presented in $\GilbargTrudinger$ to infer 
$$\sup_{X_{R-1, R+1}}|\P|\le C\|\P\|_{L^4(X_{R-2, R+2})}.$$
Combining it with (8.6) we are done. 
\qed
\enddemo

Since (8.1) is identical to (2.13) near 
time infinities, we  have 
\proclaim{Proposition 8.7} An analogue of Proposition 4.2 for (8.1) holds.
\endproclaim

We have various configuration spaces  and moduli spaces 
associated with (8.1) which are analogous to the spaces introduced in 
Sections 4 and 6. All the analysis in Sections 4, 5 and 6 
carres  over.  We shall be brief in formulating the relevant results.

Let e.g. $\R_{\pm}$ denote the $\R$ for
 $(h_{\pm}, \pi_{\pm}, \lmd_{\pm})$.
Consider 
$\a_- \in \R_{-} $, $\a_+ \in\R_{+} $ and $p_{\pm}\in\a_{\pm}$. 
We have the spaces of transition trajectories $\N(p_-, p_+)$ 
and the moduli spaces $\M^0(p_-, p_+), \M_T^0(S_{\a_-},  
S_{\a_+})$ etc..
We also have the various spaces of 
consistent multiple temporal transition  trajectory classes,
and those of consistent piecewise transition trajectories.
In particular, we have 
$\MM^0_T(S_{\a_-}, S_{\a_+})$.
For example, a multiple 
temporal transition trajectory class is a tuple 
$([u_1]_0^T,...,[u_k]_0^T) 
\in  \M^0_T(S_{\a_0^-}, S_{\a_1^-}) 
\times_{S_{\a^-_1}} ... 
\M^0_T(S_{\a^-_{m-1}}, 
S_{\a^+_{m}}) ... \M^0_T(S_{\a^+_{k-1}}, S_{\a^+_k})$  with 
a distinguished piece $[u_m]_0^T, 1 \leq m \leq k$.  
The  time translation acts on all pieces except 
the distinguished one. We define  
the $\Bbb R$-action on the distinguished one as 
the trivial action.  In particular, the $\Bbb R$-action
  on $\M^0_T(S_{\a_-}, S_{\a_+})$ is defined to be 
 the trivial action. As before, we use the underline to 
 denote quotient by the $\Bbb R^k$-action in the 
 case of temporal transition trajectories 
 and multiple temporal transition trajectory
 classes. In the case 
 of general transition trajectories and piecewise
 transition trajectories, we have the $\Bbb R^k$-action 
 induced from the twisted time translation, where the 
 $\Bbb R$-action on the  distinguished portion is 
 again trivial.

Let 
$O_{\pm}$ be the unique reducible elements 
in $\R_{\pm}$ respectively. We set 
$$
m_0= {ind~ }\F_{o_-, o_+} -1, \tag 8.9
$$
where $o_- \in O_-, o_+ \in O_+$.

The following lemma is analogous to Corollary 5.2.
\proclaim{Lemma 8.8} We have 
$$
{ind }~\F_{p_-, p_+} = \mu_-([p_-])-\mu_+([p_+])
+m_0 + {dim }~ G_{p_+}.
\tag 8.10
$$
\endproclaim

We have the following analogues 
of Theorem 4.14, Theorem 6.15 and 
Theorem 6.19.

\proclaim{Theorem 8.9} The moduli 
spaces 
${\underline \MM}^0_T(S_{\a_-}, S_{\a_+})$ 
are compact Hausdorff spaces  with respect to the 
topology induced by the 
piecewise exponential convergence. Moreover,
for generic $ b_0$, we have 
for all $\a_- \in \R_{-}, \a_+ \in \R_{+}$

(1) ${\underline {\MM}}_T^0(S_{\a_-},
S_{\a_+})$ has the structure of
$d$-dimensional smooth oriented manifolds with
corners, where 
$$d = \hbox{ind }\F_{p_-, p_+} -
max\{\hbox{dim }G_{p_-}, \hbox{dim }G_{p_+}\}+1$$
$$= \mu(\a_-)-\mu(\a_+) +m_0 + \hbox{dim }G_{p_+}
-max\{\hbox{dim }G_{p_-}, \hbox{dim }G_{p_+} \}+1$$
for $p_- \in \a_-, p_+ \in \a_+$.

(2) This structure is compatible with the natural stratification of 
${\underline {\MM}}_T^0(S_{\a_-},
S_{\a_+})$.

(3) The temporal projections $\pi_-:
{\underline {\MM}}_T^0(S_{\a_-},
S_{\a_+}) \to  S_{\a_-}$ and $\pi_+: {\underline 
{\MM}}_T^0(S_{\a_-},$
$S_{\a_+})$ $ \to S_{\a_+}$ are smooth maps. (But they 
may not be fibrations in general.)

Consequently, there holds
$$
\pa \underline \MM_T^0(S_{\a_-}, S_{\a_+}) = 
(\cup_{\mu(\a_-) > \mu(\a_-') \geq \mu(\a_+)-m_0}
\underline \MM_T^0(S_{\a_-}, S_{\a_-'}) 
\times_{S_{\a_-'}} \underline \MM_T^0(S_{\a_-'}, 
S_{\a_+}) )\cup
\tag 8.11
$$
$$ (\cup_{\mu(\a_-) \geq 
\mu(\a_+')-m_0 > \mu(\a_+)-m_0} \underline \MM_T(S_{\a_-},
S_{\a_+'}) \times_{S_{\a_+'}} \underline \MM_T^0(S_{\a_+'},
S_{\a_+})).
$$

(Of course, only the nonempty moduli spaces appear 
in this equation.)
\endproclaim

The  key point here is that Lemma 8.5 rules out 
reducible transition trajectories. Note that 
instead of using holonomy perturbations we now use 
the perturbation $b_0$ as in $\KronheimerMrowka$.
 (This perturbation 
is not time translation equivariant, and hence can't 
be applied in the construction of our homologies.)

Now we proceed to construct the desired  
shifting isomorphism.

Consider the Bott complexes $C^{Bott-}_*$ and 
$C^{Bott+}_*$ associated with the parameters
$(h_-, \pi_-, \lambda_-)$ and $(h_+, \pi_+,
\lambda_+)$ respectively.
Let   $\bold F$ denote the 
collection of all the projections 
$\pi_-: \MM_T^0(S_{\a_-}, S_{\a_+})
\to S_{\a_-}$ and 
$\pi_-: \pa \MM_T^0(S_{\a_-},
S_{\a_+}) \to S_{\a_-}$. We construct our 
chain homomorphism  $F: C^{Bott-, \bold F}_{k} \to 
 C_{k+m_0}^{Bott+}. $ 
  For simplicity, we assume that the coefficient 
group is $\Bbb Z$. The general case is similar.
Consider $\sigma\in 
 C_j^{\bold F}(S_{\a_-})$
with $\sigma = [(\Delta, f)]$ and $j+ \mu(\a_-)=k$. 
For each $\a_+ $ with the corresponding 
moduli space nonempty, we follow the construction of the 
boundary operator $\pa_{\a,\b}$ in Section 7
to obtain a generalized cubical singular  chain $\sigma' 
\in C_{j'}(S_{\a_+})$ with $j'+ \mu_+(\a_+) =
j+ \mu_-(\a_-) +m_0$.
We define $F_{\a_-, \a_+}(\sigma)$ to 
be $\sigma'$. We define it to be zero 
if the moduli space is empty. Then we set 
$$
F(\sigma) = \sum_{\mu_+(\a_+) \leq \mu_-(\a_-)+m_0}
F_{\a_-, \a_+}(\sigma).
$$

It is easy to see that  we indeed 
have $F: C_{k}^{Bott-, \bold F} \to  C_{k+m_0}^{Bott+}$. 

\proclaim{Theorem 8.10} We have 
$ \pa_{Bott} \cdot F = F \cdot  \pa_{Bott}$, hence 
$F$ is a chain map
of degree $m_0$ from $ (C_*^{Bott-, \bold F}, \pa_{Bott})$
to $(C_*^{Bott+}, \pa_{Bott})$. The induced shifting 
homomorphism between the homologies is denoted 
$F_*$. From the construction of $F$ we easily obtain 
a cochain map of degree $m_0$ and  
the induced 
shifting homomorphism $F^*$ between cohomologies. 
\endproclaim
\demo{Proof} 
The proof goes along the same lines as 
the proof of Lemma 7.5.  First, we 
arrange the orientations such that
$$
\pa \underline \MM_T^0(S_{\a_-}, S_{\a_+})
= (-1)^{\mu(\a_-)+ dim S_{\a_-}+m_0+1}
(\Sigma_{\a_-'} \underline 
\MM_T^0(S_{\a_-}, 
S_{\a_-'}) \times_{S_{\a_-'}} 
\tag 8.12
$$
$$
\underline \MM_T^0(S_{\a_-'}, S_{\a_-})  ) +
(\Sigma_{\a_+'} \underline \MM_T^0(S_{\a_-},
S_{\a_+'}) \times_{S_{\a_+'}} 
\underline \MM_T^0(S_{\a_+'}, S_{\a_+})).
$$
In comparison with (7.3), we see an additional 
$1$ appearing in the sign exponent (note that $m_0$ 
reduces to zero in the situation of (7.3)).
This is because of the difference in 
the time translation action.  
Consider 
$\sigma =
[(\Delta, f)] \in 
C_j(S_{\a_-})$.
Analogous to (7.5) we have for 
$\a_+$ with $\mu(\a_+) \leq \mu_-(\a_-)+m_0$
 
$$ \pa(\D  \times_{S_{\a_-}} \underline
{\MM}_T^0(S_{\a_-},S_{\a_+})) =
\pa
\D \times_{S_{\a_-}}\underline {\MM}_T^0(S_{\a_-},
S_{\a_+})+ 
$$
$$
(-1)^{j+dim S_{\a_-}}\D \times_{S_{\a_-}}\pa \underline
{\MM}_T^0(S_{\a_-}, S_{\a_+})
=
$$
$$
\pa \D\times_{S_{\a_-}}\underline {\MM}_{T}^0(S_{\a_-},
S_{\a_+})
+$$ 
$$ (-1)^{j+ \mu(\a_-)+ m_0 +1}\sum_{\a_-'} \D 
\times_{S_{\a_-'}}
 \underline
{\MM}_T^0(S_{\a_-},S_{\a_-'}) \times_{S_{\a_-'}}
\underline {\MM}_T(S_{\a_-'},S_{\a_+}) +$$

$$ (-1)^{j+ \mu(\a_-)+ m_0+1}\sum_{\a_+'} \D 
\times_{S_{\a_-}}
\underline
{\MM}_T^0(S_{\a_-},S_{\a_+'}) \times_{S_{\a_+'}}
\underline {\MM}_T(S_{\a_+'},S_{\a_+}).$$

This implies
$$ \pa_0 F_{\a_-, \a_+} \sigma =
F_{\a_-, \a_+} \pa_0 \sigma +\sum_{\mu(\a_-)
>\mu(\a_-') \geq \mu(\a_+)-m_0}
F_{\a_-', \a_+} \pa_{\a_-, \a_-'} 
\sigma -$$
$$
\pa_{\a_+', \a_+} F_{\a_-, \a_+'} \sigma.
$$
 Summing over all $\a_+$, we arrive 
at the desired chain map property.
\qed
\enddemo

\head Appendix A. The ordinary Seiberg-Witten Floer homology
\endhead        

Here, we only consider rational homology spheres. 
It is not hard to extend the construction to general manifolds. 

For simplicity, we formulate the theory with integer 
coefficients. Let $Y$ be a rational homology sphere 
with a given metric $h$ and $c$ a $spin^c$ structure on $Y$ as 
before. Set $C_i=\Bbb Z\{\a\in\R^*|\mu(\a)=i\}$.  We
define a boundary operator $\pa:C_i\to C_{i-1}$ in terms of the 
moduli spaces $\ti \M(\a, \b)$, or equivalently
$\ti \M(p, q), p \in \a, q\in \b$, 
where the tilde means quotient by the time translation action.
For a generic pair $(\pi, \lmd)$, $\ti \M(\a, 
\b)$ is a compact oriented manifold of zero dimension,
provided that $\mu(\a) - \mu(\b)=1$. The orientation (sign)
is given by pairing the orientation of $\M(\a, \b)$ 
provided by Proposition 5.3 with its orientation induced 
by the time translation action. For $\a \in C_i$ we then set 
$$\pa \a=\sum_{\mu(\a)-\mu(\b)=1}\sharp\ti\M(\a, \b)\b,\tag A.1$$
where $\sharp\ti\M(p,q)$ is the algebraic sum of $\ti\M(\a, \b)$.

The compactification of the moduli spaces $\ti \M(\a, \ga)$ with 
$\mu(\a) - \mu(\ga) =2$ is similar to 
the results in Section 6. For dimensional reasons, 
no trajectory connecting to 
the reducible point appears in the compactification.
Using these compactified moduli spaces and the 
consistency of orientation (Proposition 5.3) 
we obtain 

\proclaim{Lemma A.1} $\pa^2=0$.
\endproclaim

\definition{Definition A.2} The ordinary 
Seiberg-Witten Floer homology 
$FH^{SW}_*(c, h, \pi, \lmd) $ and 
cohomology $FH^{SW*}(c, h, \pi, \lmd) $ 
for the $spin^c$ structure $c$ and the parameters $h, \pi, 
\lmd$  are defined to be the homology and 
cohomology of the chain complex $(C_*, \pa)$. 
\enddefinition

\head{Appendix B. Local Gauge Fixing }
\endhead

In this appendix we prove a result on local Columb gauge fixing. 
We assume that $Y$ is a rational homology 
sphere. Let $c$ be a fixed $spin^c$ structure on $Y$ and 
$a_0$ a smooth reference connection as in Section 3. 
It induces a time-independent connection $A_0$ over 
$X$.  
For positive numbers $r < R$, let  $x_0 \in X_{r, R}$ be 
a reference point. 
 
\proclaim{Lemma B.1} For any $A\in \A_{l}(X_{r, R})$ with 
$l \geq 1$, there exists
a unique gauge $g\in \G_{l+1}(X_{r, R})$ with $g(x_0)=1$  such that
$\ti A=g^*(A)$ satisfies
$$\cases d^*(\ti A- A_0)=0,\cr
(\ti A-A_0)(\nu) =0 \hbox{ on }
\pa X_{r, R},\cr
\endcases 
\tag B.1
$$
where $\nu$ denotes the unit outer normal 
of $\pa X_{r, R}$. Moreover, we have 
$$\|\ti A -A_0\|_{1, 2}\leq C\|F_{\ti A}- F_{A_0}\|_{0, 2}$$
for a positive constant $C$ depending only on $Y, r $ and $R$.
\endproclaim
\demo{Proof}
The associated gauge fixing equation is 
$$\cases d^*(g^{-1}dg)+d^*(A-A_0) & =0,\cr
g^{-1}dg(\nu)+(A-A_0)(\nu)& =0 \hbox{ on }
\pa X_{r, R}.\cr
\endcases 
\tag B.2 
$$
If we choose $g= e^f$, then the equation reduces to 
$$
\cases d^*df+d^*(A-A_0) & =0,\cr
\frac{\pa f}{\pa \nu}+(A-A_0)(\nu)& =0 \hbox{ on }
\pa X_{r, R}.\cr
\endcases 
\tag B.3 
$$
It is clear that a solution   $f \in 
\O^0_{l+1}(X_{r, R})$ with $f(x_0) = 0$ exists.       

Now  assume that there are 
$g_1$  and $g_2$ satisfying (B.2) with 
$g_1(x_0)=g_2(x_0)=1$. We set $g = g_1^{-1} g_2$. 
Then $g$ is a solution of the problem
$$\cases d^*(g^{-1}dg) & =0,\cr
g^{-1}dg(\nu) & =0 \hbox{ on }
\pa X_{r, R}.\cr
\endcases 
$$

Taking  two copies of $X_{r, R}$ and gluing them along their 
common boundary, we obtain  the  Riemannian manifold 
$Y\times S^1$. The two copies  of $g$ then yield a solution $
g_0: Y \times S^1 \to S^1$ of the harmonic map equation $d^*(g^{-1}dg)=0$. 
This means that the 1-form $g^{-1}dg$ is harmonic.

We claim that 
every harmonic 1-form on $Y \times S^1$ is a constant 
multiple of the 
base harmonic form $ds$, where $s$ is the arclength parameter 
on $S^1$. Indeed, let $\alpha = \alpha_Y + f ds$ be an arbitary 
harmonic 
1-form, where $\alpha_Y$ has no $ds$ component.  Then we deduce 
from the equation $d \alpha =0$ 
$$
\cases d_Y \alpha_Y & =0, \cr
\frac{\partial \alpha_Y}{\partial s} & = d_Y f. \cr
\endcases
$$
Since $Y$ is a rational homology sphere, the first equation 
above implies that $\alpha_Y = d_Y h$ for a function $h$.
Then the second equation implies that $f = \frac{\partial 
h}{\partial s} + f_0$ for a function $f_0$ depending only on 
$s$.  It follows that $\alpha = d h + f_0 ds$.  Now we use 
the equation $d^* \alpha =0$ to deduce 
$$
d^* d h + \frac{d f_0}{d s}=0.
$$
Integrating $f_0$ along $S^1$ we obtain another solution 
$h_0$ of this equation, which depends only on $s$. 
Clearly
$h$ differs from $h_0$ by a constant. We deduce 
that $\frac{d h}{ ds } + f_0$ is a constant, whence 
$\alpha$ is a constant multiple of $ds$.

Since 
$g^{-1} dg$ is harmonic, it is a constant multiple of 
$ds$. On the other hand, 
for each $y\in Y$, $g_0(y,\cdot)$ is a map from $S^1$ to $S^1$
 with
degree zero, hence {\it i.e.} $\int_{\{y\}\times S^1}g_0^{-1}dg_0=0$.
it follows that $g^{-1} dg =0$.
We deduce that $g$ is constant.   Consequently,
$g_1 \equiv g_2$.

Next consider an $A$ satisfying 
(B.1).  By the above gluing argument, the form 
$A-A_0$ leads to 
a form $A_1$ on $Y\times S^1$ 
satisfying $d^*A_1=0$.  The Hodge decomposition of 
$A_1$ takes the following form
$$
A_1 = c ds + d^* \beta,
$$
where $c$ denotes a constant and $\beta \in
\Omega_2^2(Y \times S^1)$.  The constant $c$ is 
given by the $L^2$ product between $A_1$ and 
$ds$. By the construction of $A_1$, 
it is easy to see that 
this product is zero. Now we deduce
$$\|A_1\|_{1, 2}= \|d^*\beta\|_{1, 2}\le C\|d d^*\beta 
\|_{0, 2} = C\|dA_1\|_{0, 2}.$$
But $\|d A_1 \|_{0, 2} = \sqrt{2} \|F_A -F_{A_0}\|_{0,2}$,
hence the desired estimate follows.
\qed
\enddemo

\definition{Remark B.2} If $Y$ is not a homology 
sphere, then we can still achieve the gauge 
fixing (B.1) provided that the $L^2$ norm of 
the curvature of $A$ is sufficiently small. This 
is similar to Uhlenbeck's gauge fixing lemma 
in the Yang-Mills theory, cf. $\DonaldsonKronheimer$
$\Uhlenbeck$.
\enddefinition

\head{Appendix C. A transversality}
\endhead

In this appendix  we prove transversality at 
reducible trajectories. 

\proclaim{Lemma C.1} Let $(\pi, \lmd)$ be a pair 
of $Y$-generic parameters such that $\n^2 H$ is sufficiently small. 
Choose $\d_-, \d_+$ small enough (but positive)
in the set-up of Definition 4.4. Then the operator 
$\F_{p, q}$ for $p, q \in O=O_{\pi, \lambda}$
(the unique reducible element in $\R$) at  a Seiberg-Witten 
trajectory is onto.
\endproclaim
\demo{Proof} Consider $p, q \in
O$ and $\F= \F_{p, q}$ at a trajectory 
$(A_0, \P_0)$. By the arguments in the proof 
of Theorem 4.13 we 
can assume that $p=q=(a_0, 0)$ and $(A_0, \P_0)
\equiv (a_0, 0).$ 
The formal 
adjoint $\F^*$ of $\F=\F_{p, q}$ with respect to the 
product (4.4) is given by 
$$\F^*(v)=-\frac{\pa v}{\pa t}-\pmatrix
*_Ydb+d_Yf-\n^2 H(a_0)\cdot b \cr -\np_{a_0}\psi
-\lmd\psi \cr
d_Y^*b+2\d_F' f\cr
\endpmatrix-{2\d_F'}v \tag C.1
$$
for $v=(\psi,b,f)$. The surjectivity of $\F$ 
is equivalent to the vanishing of the kernel of 
$\F^*$.  Let $v$ satisfy $\F^* v=0$. 
Then we have 
$$\frac{\pa }{\pa t}\pmatrix b\cr f\cr\endpmatrix+
\pmatrix *_Y d_Y b+d_Yf\cr d_Y^*b\cr\endpmatrix+
\pmatrix {\d_F'}b+ \n^2H(a_0) \cdot b\cr 0\endpmatrix
=0. \tag C.2$$
We define the operator $L$ by 
$$L\pmatrix b\cr f\cr\endpmatrix=-\pmatrix *_Yd_Yb+d_Yf\cr d_Y^*b\cr\endpmatrix.$$
$L$ is formally self-adjoint and  satisfies 
$L^2 = \D$, where $\D$ denotes the Hodge Laplacian. 
Let
$\{\xi_i=\pmatrix b_i\cr f_i\cr\endpmatrix\}$ be a complete $L^2$ 
orthonormal system of eigenvectors of $L$ with 
$L\xi_i=\lmd_i\xi_i$. From the above discussion we deduce 
$$\D b_i=\lmd_i^2b_i,  \D f_i=\lmd_i^2 f_i.$$

Now we write 
$$\pmatrix b\cr f\cr\endpmatrix=\sum _{-\infty}^{+\infty}l_i(t)\xi_i.$$
Then it follows from C.2 that 
$$\sum l_i'(t)\xi_i+\lmd_il_i(t)\xi_i+\pmatrix 
{\d_F'}b + \n^2H(a_0) \cdot b \cr 0\endpmatrix=0.\tag C.3$$
 Consequently, $f_j$ is a nonzero 
constant. Then we deduce 
$l_{j}'(t) \equiv 0$, hence $l_j$ is a constant. But $(b, f)$ 
is $L^2$ integrable, which forces $l_j$ to be zero. We conclude 
that the above expansion of $(b ,f)$ does not contain terms 
with zero eigenvalue. Using the elementary 
arguments in e.g. $\SalamonZehnder$ it is then easy to show that $(b, f)$ 
must vanish, provided that $\n^2H$, $\d_+$ and $\d_-$ have been chosen 
small enough.  Using the same arguments one also 
infers that $\psi$ vanishes. Thus $v=0$.  
\qed
\enddemo

\bigskip
{\it Aknowledgement: The first named author was supported
by SFB 237 and the Leibniz program of DFG. He is grateful
to Prof. J. Jost for providing these supports.}


\Refs

\refstyle{1}
\widestnumber\key{JPW}
 
\ref\key 1
\by M.F. Atiyah
\paper New invariants of 3- and 4-dimensional manifolds
\jour Symp. Pure Math.
\vol 48
\yr 1988
\pages 285-299
\endref


\ref\key 2
\by M.F.Atiyah, V.K. Patodi and I.M. Singer
\paper ~Spectral asymmetry and Riemannian geometry I, II, III
\jour Math. Proc. Cambridge. Soc.
\vol 77 
\yr 
\pages 
\endref


\ref\key  3
\by D. M. Austin and P. J. Braam
\paper Equivariant Floer theory and gluing 
Donaldson polynomials 
\jour Topology
\vol 35 
\yr 1996
\pages 167-200
\endref

\ref\key 4
\bysame
\paper Morse-Bott theory and equivariant 
cohomology
\inbook Floer Memorial Volume
\bookinfo  Birkh{\"a}u-
ser 
\yr 1996
\pages 123-164
\endref

\ref\key 5
\by W. Chen
\paper  Casson's invariant and Seiberg-Witten gauge 
theory
\inbook  preprint
\bookinfo   
\yr 
\pages 
\endref

\ref\key 6
\by R. Cohen, J. Jones and G. Segal
\paper Floer's infinite dimensional Morse theory 
and homotopy theory
\inbook Floer Memorial Volume
\bookinfo Birkh{\"a}user
\yr 1996
\pages 297-326
\endref

 

\ref\key 7 
\by S.K. Donaldson
\paper The orientation of Yang-Mills moduli spaces
 and 4-manifold topology
\jour J. Diff. Geo.
\vol 26
\yr 1987
\pages 397-428
\endref


 

\ref\key 8
\bysame 
\paper  The Seiberg-Witten 
equations and the 4-manifold 
topology
\jour Bull. Amer. Math. Soc.
\vol
\yr 1996
\pages 
\endref


\ref\key 9
\by S.K. Donaldson and P. B. Kronheimer
\book The geometry of four-manifolds
\bookinfo Oxford Science Publications
\yr 1990
\endref

\ref\key 10
\by A. Floer
\paper An instanton invariant for 3-manifolds
\jour Comm. Math. Pyhs.
\vol 118
\yr 1988
\pages 215-240
\endref


\ref\key 11
\by A. Floer
\paper Holomorphic spheres and symplectic 
fixed points
\jour Comm. Math. Phys.
\vol
\yr
\pages 
\endref



\ref\key 12
\by Fukaya
\paper Floer homology of connected sum of homology 3-sphere
\jour Topology
\vol 36
\yr  1996
\pages 
\endref


\ref\key 13
\bysame 
\paper Floer homology for oriented 3-manifolds
\book Aspects of low dimensional manifolds
\bookinfo Advanced studies in pure ~mathematics
\yr 1992
\pages 1-99
\endref




\ref\key 14
\by Gilbarg and  Trudinger
\book Elliptic Partial Differential 
Equations of Second Order
\bookinfo Springer
\yr 1981
\endref

\ref\key 15
\by P. B. Kronheimer and T. S. Mrowka
\paper The genus of embedded surfaces in the projective space
\jour Math. Res. Letters
\vol 1
\yr 1994
\pages 797--808
\endref



\ref\key 16
\by  H. B. Lawson and M.-L Michelsohn
\book  Spin geometry
\bookinfo ~Princeton, New Jersey
\yr 1989
\endref

\ref\key 17
\by Y. Lim
\paper The equivalence of Seiberg-Witten and 
Casson invariants for homology 3-spheres
\jour preprint 
\yr
\endref

\ref\key 18
\by M. Marcolli 
\paper Seiberg-Witten-Floer homology and 
Heegard splittings
\jour Inter. J. Math. 
\vol 7
\yr  1996
\pages 671-696
\endref

\ref\key 19
\by W. S. Massey
\book Singular homology theory
\bookinfo Graduate Texts in Mathematics, Springer-Verlag
\yr 1980
\endref

\ref\key  20
\by D. Salamon 
\book Spin geometry and the Seiberg-Witten invariants
\bookinfo Lectures notes, Warwick
\yr 1995
\endref



\ref\key 21
\by D. Salamon and E. Zehnder
\paper Morse theory for periodic solutions 
of Hamiltonian system and
the Maslov index
\jour Comm. Pure Appl. Math.
\vol 45
\yr 1992
\pages 1303-1360
\endref


\ref\key 22 
\by C. Taubes 
\paper   Casson's invariant and gauge theory 
\jour J. Diff. Geom.
\vol 31
\yr 1990
\pages 547-599
\endref

\ref\key 23 
\bysame
\paper The Seiberg-Witten invariants and the Gromov invariants
\jour preprint 
\yr
\pages
\endref




\ref\key 24
\by K. K. Uhlenbeck
\paper Connections with $L^p$ bounds on curvature
\jour Comm. Math. Phys. 
\vol  83
\yr 1982
\pages 31-42
\endref
 
\ref\key 25
\by R. Wang
\paper On Seiberg-Witten Floer invariants and 
the generalized Thom problem
\jour preprint
\vol
\yr 
\pages 
\endref


\ref\key 26
\by G. Wang and Rugang Ye
\paper Bott-type and equivariant Seiberg-Witten Floer 
homology I
\jour DG9701010
\vol 
\yr 1997
\pages 
\endref

\ref\key 27
\by E. Witten
\paper Monopoles and 4-manifolds
\jour Math. Res. Letters
\vol 1
\yr 1994
\pages 769--796
\endref

\ref\key 28
\by Rugang Ye 
\paper Equivariant and Bott-type Seiberg-Witten 
Floer homology: Part II
\jour preprint 
\vol 
\yr  
\pages 
\endref

\ref\key 29
\by Rugang Ye 
\paper Equivariant and Bott-type Seiberg-Witten 
Floer homology: Part III
\jour in preparation 
\vol 
\yr  
\pages 
\endref
 



\endRefs
\enddocument